\documentclass[twoside, 12pt]{paper}
\usepackage{a4wide, amsthm, amssymb, amscd, concmath, euler, eucal}

\DeclareSymbolFont{rsfs}{OMS}{rsfs}{m}{n}
\DeclareSymbolFontAlphabet{\mathscr}{rsfs}
\DeclareSymbolFont{bbold}{U}{dsrom}{m}{n}
\DeclareSymbolFontAlphabet{\mathbb}{bbold}

\renewcommand{\rho}{\varrho}
\newcommand{\A}{{\cal A}}

\renewcommand{\Bbb}{\mathbb}
\renewcommand{\frak}{\mathfrak}
\newcommand{\catqot}{/\hskip-3pt/}
\newcommand{\C}{{\Bbb C}}

\newcommand{\E}{{\cal E}}

\newcommand{\GL}{\mathop{\rm GL}}

\newcommand{\Hom}{\mathop{\rm Hom}}

\newcommand{\id}{\mathop{\rm id}}

\renewcommand{\L}{{\cal L}}

\newcommand{\n}{{\cal N}}
\renewcommand{\O}{{\cal O}}

\renewcommand{\P}{{\Bbb P}}

\newcommand{\Q}{{\Bbb Q}}
\newcommand{\R}{{\Bbb R}}

\newcommand{\SL}{\mathop{\rm SL}}

\newcommand{\Spec}{\mathop{\rm Spec}}

\newcommand{\Sym}{\mathop{\rm Sym}}

\renewcommand{\tilde}{\widetilde}

\newcommand{\Z}{{\Bbb Z}}

\newcommand{\la}{\lambda}
\newcommand{\lra}{\longrightarrow}
\newcommand{\ra}{\rightarrow}

\newcommand{\lma}{\longmapsto}
\newcommand{\p}{\prime}
\newcommand{\q}{\quad}

\renewcommand{\phi}{\varphi}
\newcommand{\rk}{\mathop{\rm rk}}
\newcommand{\eps}{\varepsilon}

\renewcommand{\theta}{\vartheta}
\newcommand{\ul}{\underline}
\newcommand{\ol}{\overline}

\theoremstyle{plain}
\newtheorem{Thm}{\sc Theorem}[subsection]
\newtheorem*{Thm*}{\sc Main Theorem}

\newtheorem{Prop}[Thm]{\sc Proposition}
\newtheorem*{Prop*}{\sc Proposition}
\newtheorem{Lem}[Thm]{\sc Lemma}

\newtheorem*{Claim}{\sc Claim}
\newtheorem*{Conc}{\sc Conclusion}

\theoremstyle{definition}

\theoremstyle{remark}

\newtheorem{Rem}[Thm]{Remark}
\newtheorem{Ex}[Thm]{Example}
\newtheorem*{Rem*}{Remark}

\begin{document}

\pagestyle{myheadings}
\markboth{\rm Alexander H.W.\ Schmitt}{\rm Semistable Singular Principal Bundles}

\title{A closer look at semistability for singular principal bundles\footnote{To appear in the International
Mathematics Research Notices}}
\author{Alexander H.W.\ Schmitt
\institution{\rm \small Universit\"at Duisburg-Essen\\\rm \small FB6 Mathematik \& Informatik\\ \rm\small D-45117 Essen\\ \rm\small Germany
\\ \small\tt alexander.schmitt@uni-essen.de\rm}}
\date{}
\maketitle
\begin{abstract}
We substantially refine the theory of singular principal bundles introduced in a former paper. In particular,
we show that we need only honest singular principal bundles in our compactification. These are objects which
carry the structure of a rational principal bundle in the sense of Ramanathan. Moreover, we
arrive at a much simpler semistability condition. In the case of a semisimple group, this is just the Gieseker-version
of Ramanathan's semistability condition for the corresponding rational principal $G$-bundle.
\end{abstract}
\section{Introduction}
In our paper \cite{Schmitt}, we presented an approach for compactifying the moduli spaces of semistable
principal $G$-bundles over a polarized higher dimensional base manifold $(X,\allowbreak\O_X(1))$, 
for $G$ a reductive linear algebraic group. 
For this, we fixed a faithful representation
$\rho\colon G\lra \GL(V)$ with $\rho(G)\subset \SL(V)$. Then, we looked at pairs $({\cal A},\tau)$
with ${\cal A}$ a torsion free sheaf of rank $\dim(V)$ with trivial determinant and $\tau\colon
\Sym^*(V\otimes {\cal A})^G\lra \O_X$ a homomorphism of $\O_X$-algebras which is non-trivial in the sense that the
induced section $\sigma\colon X\lra \ul{\Spec}(\Sym^*(V\otimes {\cal A})^G)$ be not the zero section. 
Such a pair was called a \it singular principal $G$-bundle\rm, and if, furthermore, $\sigma_U(U)\subset
\ul{\rm Isom}(V\otimes\O_U,{\cal A}^\vee_{|U})$, we spoke of an \it honest singular principal $G$-bundle\rm.
Here, $U$ is the maximal open subset over which ${\cal A}$ is locally free. In the case of an honest
singular principal $G$-bundle $({\cal A},\tau)$, we get a principal $G$-bundle $\cal P({\cal A},\tau)$
over $U$, defined by means of base change:
$$
\begin{CD}
\cal P({\cal A},\tau) @>>> \ul{\rm Isom}(V\otimes\O_U,{\cal A}^\vee_{|U})
\\
@VVV @VVV
\\
U @>\sigma_U >> \ul{\rm Isom}(V\otimes\O_U,{\cal A}^\vee_{|U})/G.
\end{CD}
$$
For any positive polynomial $\delta\in\Q[x]$ of degree at most $\dim (X)-1$, we obtained the notion 
of \it $\delta$-(semi)stability \rm for singular principal $G$-bundles, and we managed to construct projective
moduli spaces for $\delta$-semistable principal $G$-bundles $({\cal A},\tau)$ with fixed Hilbert polynomial 
$P=P({\cal A})$.
These definitions and results have two main drawbacks:
\begin{itemize}
\item It is not clear that we need only honest singular principal bundles in order to obtain a projective
      moduli space.
\item The notion of $\delta$-semistability seems complicated and unnatural.
\end{itemize}
Using results of G\'omez and Sols (\cite{GS}, \cite{GS2}), one can establish the following
properties for $G$ a classical group and $\rho$ its ``standard" representation or $G$ an adjoint group
and $\rho$ the adjoint representation, and $\delta$ a polynomial of degree exactly $\dim(X)-1$:
\begin{itemize}
\item One needs only honest singular principal bundles in order to obtain a projective moduli space.
\item The concept of $\delta$-semistability does not depend on $\delta$, it implies the Mum\-ford-semistability
      of ${\cal A}$, and is a natural generalization of Rama\-nathan's notion of semistability.
\end{itemize}
The strategy of G\'omez and Sols to derive these properties is the following :
There is a representation $\kappa\colon \GL(V)\lra \GL(W)$ for which
there is a point $w_0\in W$, such that  $G$ has finite index in $\widetilde{G}$, the stabilizer of $w_0$, 
and such that, for every one parameter subgroup $\la\colon \C^*\lra \GL(W)$ one has $\mu(\la, w_0)\ge 0$
with equality if and only if $\la$ is actually conjugate to a one parameter subgroup of $G$
within the associated parabolic subgroup $Q_{\GL(V)}(\la)$, see~(\ref{parabolic}).
For example, if $G={\rm SO}(r)$, we have $V=\C^r$ and $W=\{\,$symmetric $(r\times r)$-matrices$\,\}$,
$w_0={\Bbb E}_r$, and $\widetilde{G}={\rm O}(r)$, or, if $G$ is an adjoint group, then $V={\frak g}$ and 
$W=\Hom(\frak g\otimes\frak g,
\frak g)$, $w_0=\hbox{Lie-bracket of $\frak g$}$, and $\widetilde{G}={\rm Aut}(\frak g)$. 
In this paper, we will extend the scope of this strategy by showing that we can, in fact, 
always find a representation $\kappa$ with the necessary conditions to establish the above properties.
Since the representation $\kappa$ will be given only abstractly, the methods we use will, however, be very different
from those of G\'omez and Sols.
\par
We will now describe the resulting notion of semistability
in more detail. Let $({\cal A},\tau)$ be an honest singular
principal $G$-bundle, and $\la\colon \C^*\lra G$ a one parameter subgroup of $G$. Recall that this yields
a parabolic subgroup $Q_G(\la)$ (see~(\ref{parabolic}) below) and a weighted flag $(V^\bullet, \ul{\alpha})$
in $V$ (see~\ref{weighted}, iv). Then, a \it reduction of $({\cal A},\tau)$ to $\la$ \rm is a section
$\beta\colon U^\p \lra {\cal P}({\cal A},\tau)_{|U^\p}/Q_G(\la)$ over an open subset $U^\p\subset U$ with
${\rm codim}_X(X\setminus U^\p)\ge 2$. This defines a weighted filtration $({\cal A}^\bullet_\beta,\ul{\alpha}_\beta)$
of $\A$.
Here, $\ul{\alpha}_\beta=(\alpha_{s},...,\alpha_1)$, if $\ul{\alpha}=(\alpha_1,...,\alpha_s)$, and the filtration 
${\cal A}_\beta^\bullet\colon\ 0\subsetneq {\cal A}_1\subsetneq \cdots\subsetneq {\cal A}_s\subsetneq {\cal A}$ is obtained as 
follows:
The section
$$
\beta^\p\colon\ U^\p\stackrel{\beta}{\lra} {\cal P}({\cal A},\tau)_{|U^\p}/Q_G(\la)\hookrightarrow
\ul{\rm Isom}(V\otimes\O_{U^\p},{\cal A}_{|U^\p}^\vee)/Q_{\GL(V)}(\la)
$$   
yields a filtration 
$$
0\subsetneq {\cal A}^\p_1\subsetneq \cdots \subsetneq {\cal A}^\p_s\subsetneq {\cal A}^\vee_{|U^\p}
$$
of ${\cal A}^\vee_{|U^\p}$ by subbundles with $\rk(\A_i^\p)=\dim(V_i)$, $i=1,...,s$. This is because
$Q_{\GL(V)}(\la)$ is just the $\GL(V)$-stabilizer of the flag $V^\bullet$ and, thus,
$\ul{\rm Isom}(V\otimes\O_{U^\p},{\cal A}_{|U^\p}^\vee)/Q_{\GL(V)}(\la)\lra U$ is the bundle
of flags in the fibres of $\A_{|U^\p}^\vee$ having the same dimensions as the flag $V^\bullet$.
We define $\A^{\p\p}_i:=\ker(\A_{|U^\p}\lra \A_{s+1-i}^{\p^\vee})$, $i=1,...,s$,
so that we obtain
a filtration 
$$
0\subsetneq {\cal A}^{\p\p}_1\subsetneq \cdots \subsetneq {\cal A}^{\p\p}_s\subsetneq {\cal A}_{|U^\p}
$$
of ${\cal A}_{|U^\p}$ by subbundles. Note that
\begin{equation}
\label{degreeequality}
\deg(\A^{\p\p}_i) = \deg({\A}^\p_{s+1-i}),\q i=1,...,s,
\end{equation}
w.r.t.\ any polarization of $X$, because $\det(\A)\cong\O_X$.
Let $j\colon U^\p\lra X$ be the inclusion and define
${\cal A}_i$ as the saturation of ${\cal A}\cap j_*({\cal A}^{\p\p}_i)$, $i=1,...,s$.
It is worth noting that, if $\la^\p=g\cdot\la\cdot g^{-1}$ for some $g\in G$, then any reduction to $\la$
may also be interpreted as a reduction to $\la^\p$.
Now, we say that an honest singular principal $G$-bundle $({\cal A},\tau)$ is \it (semi)stable\rm,
if for every one parameter subgroup $\la\colon \C^*\lra G$ and every reduction $\beta$ of $({\cal A},\tau)$
to $\la$, we have
	  $$
	  M({\cal A}^\bullet_\beta, \ul{\alpha}_\beta)\q (\succeq)\q 0.
	  $$
Recall from \cite{Schmitt} that, for every weighted filtration $({\cal A}^\bullet, \ul{\alpha})$
of ${\cal A}$,
$$
M({\cal A}^\bullet, \ul{\alpha}):=\sum_{i=1}^s\alpha_i\bigl(P({\cal A})\rk \A_i-P({\cal A}_i)\rk\A\bigr).
$$
\begin{Rem*}
i) The stated condition implies that the sheaf $\A$ is Mumford-semistable.
\par
ii) If $G$ is semisimple, then we have the implications
$$
\begin{array}{rcl}
\cal P({\cal A},\tau) \hbox{ is Ramanathan-stable} &\Longrightarrow& ({\cal A},\tau) \hbox{ is stable}
\\
&\Longrightarrow& ({\cal A},\tau) \hbox{ is semistable}
\\
&\Longrightarrow& \cal P({\cal A},\tau) \hbox{ is Ramanathan-semistable}.
\end{array}
$$
More precisely, in our language, Ramanathan's notion of (semi)stability becomes  
\begin{equation}
\label{slOpe}
\sum_{i=1}^s\alpha_i\bigl(\deg({\cal A})\rk \A_i-\deg({\cal A}_i)\rk\A\bigr)(\ge)0
\end{equation}
for every one parameter subgroup $\la\colon \C^*\lra G$ and every reduction $\beta$ of $({\cal A},\tau)$
to $\la$. Here, $\deg$ stands for the degree w.r.t.\ the chosen polarization.
Thus, in our (semi)stability concept we have just replaced degrees by Hilbert polynomials
whence our (semi)stability concept might be viewed as a reasonable Gieseker-version of Rama\-nathan-(semi)stability.
The detailed discussion is contained in Section~\ref{AnStab}.
\par
iii) For reductive groups other than semisimple ones, our notion of slope-semistability is more restrictive than 
Ramanathan's. In fact, we require Equation (\ref{slOpe}) to hold for reductions to \sl any \rm one parameter
subgroup whereas Ramanathan looks only at reductions to one parameter subgroups of the commutator subgroup
$[G,G]$.
The difference might be best understood for the reductive group $G:=\prod_{i=1}^t\GL_{r_i}(\C)$. 
Then, we may identify a $G$-bundle with a tuple of vector bundles $(E_i, i=1,...,t)$. 
Ramanathan's concept of (semi)stability just says that each $E_i$ is a (semi)stable vector bundle, $i=1,...,t$.
Our notion of slope-semistability is equivalent to the fact that $E_1\oplus\cdots\oplus E_t$ is slope-semistable
(by i)).
Therefore, it might happen that we find no semistable objects although there are Ramanathan semistable objects.
Thus, in that case, our results will be only for some special topological invariants an alternative to the
work of G\'omez and Sols.
\end{Rem*}
As usual, we define moduli functors
$$
\begin{array}{rrcl}
\ul{\hbox{M}}(\rho)_{P}^{\rm (s)s}\colon & \ul{\hbox{Sch}}_\C &\lra & \ul{\hbox{Set}}\\
& S & \lma & \left\{
\begin{array}{l}
\hbox{Equivalence classes of families of}
\\
\hbox{(semi)stable honest singular principal}
\\
\hbox{$G$-bundles with Hilbert polynomial $P$}
\end{array}\right\}.
\end{array}
$$
We then have
\begin{Thm*}
There exist a projective scheme ${\cal M}(\rho)_P^{\rm ss}$ and an open subscheme
${\cal M}(\rho)_P^{\rm s}\subset {\cal M}(\rho)_P^{\rm ss}$ as well as natural
transformations of functors
$$
\theta^{\rm(s)s}\colon\ \ul{\hbox{\rm M}}(\rho)_P^{\rm (s)s}\lra h_{{\cal M}(\rho)_P^{\rm (s)s}}
$$
with the following properties:
\begin{enumerate}
\item[\rm 1.] For every scheme $\n$ and every natural transformation $\theta^\p
\colon {\cal M}(\rho)_P^{\rm ss} \lra h_\n$, there is one and only one morphism
$\psi\colon {\cal M}(\rho)_P^{\rm ss}\lra \n$ with
$\theta^\p= h(\psi)\circ \theta^{\rm ss}$.
\item[\rm 2.] The scheme ${\cal M}(\rho)_P^{\rm s}$ is a coarse moduli space for the functor
$\ul{\hbox{\rm M}}(\rho)_{P}^{\rm s}$.
\end{enumerate}
\end{Thm*}
\begin{Rem*}
i) Note that, if $G$ is a group of the adjoint type and $\kappa$ is the adjoint representation, 
this result is contained in the work of G\'omez and Sols.
\par
ii) Balaji has recently established criteria for the non-emptiness of the
${\cal M}(\rho)_{P}^{\rm s}$ when $X$ is a surface. 
We refer to his forthcoming paper \cite{Bal}.
\end{Rem*}
We hope that the results of this paper will make the theory developed in \cite{Schmitt} 
more transparent and applicable.   Finally, let us mention
that our approach has the following advantages:
\begin{itemize}
\item We get an alternative construction for the moduli space of principal $G$-bundles over curves,
      in case $G$ is a semisimple group
      (compare  \cite{Ramanathan}, \cite{Faltings}, and \cite{Balaji}). 
	  In fact, the techniques used here may be viewed as an alternative approach
	  to the semistable reduction theorem (\cite{Faltings}, \cite{Balaji}) which readily extends to
	  the semistable reduction theorem for semistable honest singular $G$-bundles in higher dimensions.
	  On the other hand, the approach of Balaji and Seshadri yields the semistable reduction theorem
	  for \sl slope\rm-semistable honest singular $G$-bundles in higher dimensions. When $X$ is a surface,
	  the latter result enables one to give an algebraic construction of the Donaldson-Uhlenbeck compactification.
	  This is all explained in \cite{Bal}. 
\item It might be also applied to singular varieties \cite{Bhosle}, if one can make sense of the condition
      ``$\det \A$ is trivial". For recent progress in the case of irreducible nodal curves, 
	  we refer the reader to \cite{Sun}. In our forthcoming paper \cite{SchGPB}, we will use the results
	  of this and our former paper \cite{Schmitt} to obtain 
	  ``nice'' moduli spaces for singular principal bundles on irreducible nodal curves.
\item Using a faithful representation $\rho\colon G\lra \GL(V)$ allows one to treat ``decorated" singular principal
      bundles as well. This is because any representation of $G$ extends to a representation of $\GL(V)$
	  (see \cite{Deligne}, p.~40). We intend to treat this theory (over curves) in the future.
\end{itemize}
\section*{Conventions}
The general setting is as in \cite{Schmitt}.
We work over the field of complex numbers.
A \it scheme \rm will be a scheme of finite type over $\C$.
For a vector bundle $\E$ over a scheme $X$, we set $\P(\E):={\rm Proj}(\Sym^*(\E))$, i.e., $\P(\E)$ is the projective
bundle of hyperplanes in the fibres of $\E$. An open subset $U\subset X$ is said to be \it big\rm, if
${\rm codim}_X(X\setminus U)\ge 2$.
\section*{Acknowledgment}
The author's special thanks go to H.\ Kraft and J.\ Kuttler for providing the argument
for Example \ref{HomSpace} which filled a gap in an earlier version of the paper.
The author acknowledges support by the DFG through a Heisenberg fellowship and through the ``Schwerpunkt" program
``Globale Methoden in der Komplexen Geometrie --- Global Methods in Complex Geometry''.
Final revisions were made during the author's visit to the Consejo Superior de Investigaciones Cient\'\i ficas (CSIC) in Madrid
which was funded by the European Differential Geometry Endeavour (EDGE),
EC FP5 contract no. HPRN-CT-2000-00101. 
The author wishes to thank O.\ Garc\'\i a-Prada for the invitation and hospitality.
Last but not least, the referee made
useful remarks for improving the exposition.
\section{Preliminaries}
\subsection{Geometric Invariant Theory}
\label{GIT}
Let $G$ be a complex reductive group which acts on the projective scheme $X$, and suppose this action is linearized
in the ample line bundle $\L$. Given a one parameter subgroup $\la\colon \C^*\lra G$ and a point $x\in X$,
we form $x_\infty:=\lim_{z\ra \infty}\la(z)\cdot x$. Then, $x_\infty$ remains fixed under the $\C^*$-action
induced by $\la$ and the $G$-action, so that $\C^*$ acts on $\L\langle x\rangle$ by a character, say, $z\lma z^\gamma$,
$z\in\C^*$. One sets
$$
\mu_\L(\la,x)\q:=\q -\gamma.
$$ 
\begin{Lem}
\label{equality500}
Let $G$ be a reductive algebraic group, $X$ and $Y$ projective schemes equipped with
a $G$-action, and $\pi\colon X\lra Y$ a finite and $G$-equivariant morphism. Suppose
$\L$ is a $G$-linearized ample line bundle on $Y$. Then, for any point $x\in X$ and any
one parameter subgroup $\la\colon\C^*\lra G$, one has
$$
\mu_{\pi^*\L}(\la,x)\q=\q \mu_\L(\la,\pi(x)).
$$
\end{Lem}
\begin{proof}
Without loss of generality, we may suppose that $\L$ and $\pi^*\L$ are both very ample.
Define $V:= H^0(X,\pi^*\L)$ and $W:=H^0(Y,\L)$. These are $G$-modules, the inclusion
$\iota\colon W\subset V$ is $G$-equivariant and yields a $G$-equivariant rational map 
$\ol{\pi}\colon \P(V)\dasharrow \P(W)$,
and there is the following commutative diagram
$$
\begin{CD}
 X @>\hookrightarrow >> \P(V)
 \\
 @V \pi VV  @VV \ol{\pi} V
 \\
 Y @>\hookrightarrow >> \P(W)
\end{CD}
$$
of $G$-equivariant maps. Now, choose a $G$-module splitting $V^\vee \cong W^\vee \oplus \ker(\iota^\vee)$.
Note that $X\cap \P\bigl(\ker(\iota^\vee)^\vee\bigr)=\varnothing$, because $\ol{\pi}$ is defined in $X$.
A one parameter subgroup $\la$ defines splittings into non-trivial eigenspaces
$$
W^\vee \cong \bigoplus_{i=1}^l W_i\q\hbox{and}\q \ker(\iota^\vee)\cong \bigoplus_{j=1}^m I_j.
$$
For ${\Bbb V}=W_1,...,W_l, I_1,...,I_m$, let $\gamma({\Bbb V})$ be the integer with
$\la(z)\cdot [v]=[z^{\gamma({\Bbb V})}\cdot v]$ for all $z\in\C^*$, $v\in {\Bbb V}\setminus \{0\}$.
Then, for $x=[v]$,
\begin{eqnarray*}
\mu_{\pi^*\L}(\la, x) &=& \max\bigl\{\,\gamma({\Bbb V})\,|\, v \hbox{ has a non-trivial component in
${\Bbb V}$}\,\bigr\} 
\\
\mu_{\L}(\la, \pi(x)) &=& \max_{{\Bbb V}=I_1,...,I_m}
\bigl\{\,\gamma({\Bbb V})\,|\, v \hbox{ has a non-trivial component in
${\Bbb V}$}\,\bigr\} .
\end{eqnarray*}
Therefore, $\mu_{\L}(\la, \pi(x))\le \mu_{\pi^*\L}(\la, x)$.
Suppose $\mu_{\pi^*\L}(\la,x)=\gamma$. Then, there must be an index $1\le i_0\le l$, such that
$W_{i_0}$ is the eigenspace for the character $z\lma z^\gamma$ and $v$ has a non-trivial projection to
$W_{i_0}$, where $[v]=x$. Otherwise, we would have
$$
\lim_{z\rightarrow \infty} \la(z)\cdot x\in X\cap \P\bigl(\ker(\iota^\vee)^\vee\bigr),
$$
but the right hand side is empty. This shows that $\mu_\L(\la,\pi(x))\ge \gamma$ and finishes the proof.
\end{proof}
\paragraph{One Parameter Subgroups and Parabolic Subgroups. ---}
Let $G$ be a complex reductive group, and $\la\colon \C^*\lra G$ a one parameter subgroup.
Then, we define the parabolic subgroup
\begin{equation}
\label{parabolic}
Q_G(\la)\q :=\q \Bigl\{\, g\in G\,|\, \lim_{z\rightarrow\infty} \la(z)\cdot g\cdot\la(z)^{-1} \hbox{ exists in $G$}
\,\Bigr\}.
\end{equation}
In fact, any parabolic subgroup of $G$ arises in this way. We refer the reader to the books \cite{Spri} and
\cite{GIT}, Chapter 2.2, for more details. The centralizer $L_G(\la)$ of $\la$ is a \it Levi-component 
of $Q_G(\la)$\rm, i.e., $Q_G(\la)={\cal R}_u(Q_G(\la))\rtimes L_G(\la)$. In this picture, the unipotent radical
of $Q_G(\la)$ is characterized as
$$
{\cal R}_u(Q_G(\la))\q =\q \Bigl\{\, g\in Q_G(\la)\,|\, 
\lim_{z\rightarrow\infty} \la(z)\cdot g\cdot\la(z)^{-1} =e
\,\Bigr\}.
$$
\begin{Rem}
\label{weighted}
i) In the sources quoted above, one takes the limit $z\rightarrow 0$ in order to define a parabolic subgroup
$P_G(\la)$. Thus, we have
\begin{equation}
\label{Sprcomparison}
Q_G(\la)\q=\q P_G(-\la).
\end{equation}
\par
ii) Let $G$ be a complex reductive group which acts on the projective scheme $X$, and suppose this action is linearized
in the ample line bundle $\L$. Then, for any point $x\in X$, any one parameter subgroup $\la\colon\C^*\lra G$,
and any $g\in Q_G(\la)$
$$
\mu_\L(\la, x)\q=\q \mu_\L(\la,g\cdot x).
$$
This is proved in \cite{GIT}, Chapter 2.2.
\par
iii)
If we are given an injective homomorphism $\iota\colon G\hookrightarrow H$, then we obviously find
$$
Q_H(\la)\cap G\q=\q Q_G(\la).
$$
\par
iv) If $G=\GL(V)$, then the group $Q_G(\la)$ is the stabilizer of the flag
$$
V^\bullet\colon\ 0\subsetneq V_1\subsetneq V_2\subsetneq\cdots\subsetneq V_s\subsetneq V 
$$
where $V_i:=\bigoplus_{j=1}^i V^j$, $V^j$ is the eigenspace of the $\C^*$-action coming from $\la$ for
the character $z\lma z^{\gamma_j}$, and $\gamma_1<\cdots<\gamma_{s+1}$ are the different weights occurring.
We also set $\alpha_i:=(\gamma_{i+1}-\gamma_i)/\dim(V)$, $i=1,...,s$. The pair $(V^\bullet,\ul{\alpha})$
is referred to as the \it weighted flag of $\la$\rm. Note, that if $\la^\p$ is conjugate to $\la$, then
$\dim(V_i^\p)=\dim(V_i)$ and $\alpha_i^\p=\alpha_i$, $i=1,...,s$.
\end{Rem}
\begin{Ex}[Actions on Homogeneous Spaces]
\label{HomSpace}
Let $H$ be a reductive algebraic group, $G$ a closed reductive subgroup, and $X:=H/G$
 the associated affine homogeneous space. Then, the following holds true:
\begin{Prop}
Suppose that we are given a point $x\in X$ and a one parameter subgroup
$\la\colon \C^*\lra H$, such that $x_0:=\lim_{z\ra\infty} \la(z)\cdot x$ exists in $X$.
Then, $x\in {\cal R}_u(Q_{H}(\la))\cdot x_0$.
\end{Prop}
\begin{proof}[Proof (after Kraft/Kuttler)]
We may assume $x_0=[e]$. Define 
$$
Y:=\bigl\{\,y\in X\,|\,\lim_{z\ra \infty}\la(z)\cdot y=x_0\,\bigr\}.
$$
This set is closed and invariant under the action of ${\cal R}_u(Q_{H}(\la))$.
Note that viewing $X$ as a variety with $\C^*$-action, $x_0$ is the unique point in $Y$ with a closed
$\C^*$-orbit, and by the first lemma in Section III of \cite{Luna}, there
is a $\C^*$-equivariant morphism $f\colon X\lra T_{x_0}(X)$ which maps $x_0$ to $0$
and is \'etale in $x_0$. Obviously, $f$ maps $Y$ to
\begin{equation}
\label{TangComp}
\bigl\{\, v\in T_{x_0}X\,|\, \lim_{z\ra \infty} \la(z)\cdot v=0\,\bigr\}
={\frak u}_H(\la)/{\frak u}_G(\la)\subset {\frak h}/{\frak g}.
\end{equation}
Here, ${\frak u}_H(\la)$ and ${\frak u}_G(\la)$ are the Lie algebras of the unipotent radicals
${\cal R}_u(Q_H(\la))$ and ${\cal R}_u(Q_G(\la))$, respectively, and ${\frak h}$ and ${\frak g}$
are the Lie algebras of $H$ and $G$, respectively. 
Note that ${\frak h}$ and ${\frak g}$ receive their $G$-module structures through
the adjoint representation of $G$, and, moreover, by definition,
$$
{\frak u}_H(\la)=\bigl\{\,v\in {\frak h}\,|\,\lim_{z\ra \infty} \la(z)\cdot v=0\,\bigr\}.
$$
This yields the asserted equality in (\ref{TangComp}). 
The morphism $f$ provides a $\C^*$-equivariant isomorphism $Y\lra {\frak u}_H(\la)/{\frak u}_G(\la)$.
(e.g., Theorem 3.4 in \cite{Hesselink}).
On the other hand, ${\frak u}_H(\la)/{\frak u}_G(\la)$ equals the tangent space
of the ${\cal R}_u(Q_H(\la))$-orbit of $x_0$ at $X$. Therefore, since
${\cal R}_u(Q_H(\la))\cdot x_0\subset Y$, $Y$ must agree with the closed orbit 
${\cal R}_u(Q_H(\la))\cdot x_0$, and we are done.
\end{proof}
\end{Ex}
\paragraph{The Instability Flag. ---}
In this section, ${\Bbb K}$ will be an algebraically closed field of characteristic zero. (Besides for $\C$,
we will need the results also for the algebraic closure of the function field of $X$.)
We start with the group $\GL_n({\Bbb K})$. Let $T$ be the maximal torus of diagonal matrices.
The characters $e_i\colon {\rm diag}(l_1,...,l_n)\lma l_i$, $i=1,...,n$, form a basis for the character
group $X^*(T)$, and 
$$
\begin{array}{rccc}
(.,.)^*\colon & X^*_\R(T)\times X^*_\R(T) &\lra& \R
\\
& \Bigl(\sum_{i=1}^n x_i\cdot e_i, \sum_{i=1}^n y_i\cdot e_i\Bigr) &\lma &\sum_{i=1}^n x_iy_i
\end{array}
$$
defines a scalar product on $X^*_\R(T):=X^*(T)\otimes_\Z\R$ which is invariant under the action of the Weyl group
$W(T):={\cal N}(T)/T$. This yields isomorphisms
$$
X^*_\R(T)\q\cong\q {\Hom}_\R(X^*_\R(T),\R)\q\cong\q X_{*,\R}(T)\q:=\q X_*(T)\otimes_\Z\R.
$$
For the second identification, we use the duality pairing $\langle.,.\rangle_\R\colon X_{*,\R}(T) \times X^*_\R(T)
\lra \R$ which is the $\R$-linear extension of the canonical pairing 
$\langle.,.\rangle\colon X_{*}(T) \times X^*(T)
\lra \Z$. Since the pairing $(.,.)^*$ is $W(T)$-invariant, the norm $\|.\|_*$ induced on $X_{*,\R}(T)$ extends
to a $\GL_n({\Bbb K})$-invariant norm $\|.\|$ on the set of all one parameter subgroups of $\GL_n({\Bbb K})$ (see \cite{GIT},
Chapter 2.2, Lemma 2.8).
\par
Next, suppose we are given a representation $\kappa\colon \GL_n({\Bbb K})\lra \GL(W)$. This leads to a decomposition
$$
W\q\cong\q \bigoplus_{\chi\in X^*(T)} W^\chi
$$
of $W$ into eigenspaces and defines the \it set of weights of $\kappa$ (w.r.t.\ $T$)\rm
$$
{\rm WT}(\kappa,T)\q:=\q\bigl\{\,\chi\in X^*(T)\,|\, W^\chi\neq\{0\}\,\bigr\},
$$
and, for any $w\in W$, the \it set of weights of $w$ (w.r.t.\ $T$)\rm
$$
{\rm WT}(w,T)\q:=\q\bigl\{\,\chi\in {\rm WT}(\kappa,T)\,|\, w \hbox{ has a non-trivial component in
}W^\chi\,\bigr\}.
$$
For a one parameter subgroup $\la\in X^*_\R(T)$, we then set
$$
\mu_\kappa(\la,w)\q:=\q \max\bigl\{\,\langle\la,\chi\rangle_\R\,|\, \chi\in {\rm WT}(w,T)\,\bigr\}.
$$
For any other maximal torus $T^\p\subset G$, we choose an element $g\in G$ with $g\cdot T^\p\cdot g^{-1}=T$, and
set, for $\la\in X^*_\R(T^\p)$,
\begin{equation}
\label{DeFine}
\mu_\kappa(\la,w)\q:=\q \mu_\kappa(g\cdot\la\cdot g^{-1}, g\cdot w).
\end{equation}
\begin{Ex}
i) Let $\P(W^\vee)$ denote the space of lines in $W$. Then, $\kappa$ yields an action of $\GL_n({\Bbb K})$ on 
$\P(W^\vee)$
and a linearization of that action in $\O_{\P(W^\vee)}(1)$. With the former notation, we find
$$
\mu_{\O_{\P(W^\vee)}(1)}(\la, [w])\q=\q\mu_\kappa(\la, w),
$$
for every point $w\in W\setminus\{0\}$ and every one parameter subgroup 
$\la\colon {\Bbb G}_m({\Bbb K})\lra \GL_n({\Bbb K})$.
\par
ii) Our convention is the same as in \cite{Schmitt0} and \cite{Schmitt}, but differs from the one in
\cite{RamRam}. More precisely, let $\mu^{\rm RR}_\kappa(\la, w)$ be the quantity defined in \cite{RamRam}.
Then,
\begin{equation}
\label{RRcomparison}
\mu_\kappa(\la, w)\q=\q -\mu_\kappa^{\rm RR}(-\la,w).
\end{equation}
\end{Ex}
Now, suppose we are also given a reductive subgroup $G\subset \SL_n({\Bbb K})$. For simplicity, assume that there is a 
maximal torus $T_G$ of $G$ which is contained in  $T$. Otherwise, we may pass to a different maximal torus 
$T^\p$ of $\GL_n({\Bbb K})$.
From $(.,.)^*$ and the dual pairing $(.,.)_*\colon\ X_{*,\R}(T)\times X_{*,\R}(T)\lra \R$, we obtain the
induced pairing $(.,.)_{*,G}\colon\ X_{*,\R}(T_G)\times X_{*,\R}(T_G)\lra \R$. Let $\|.\|_G$ be the restriction of the
norm $\|.\|$ to the one parameter subgroups of $G$. Note that, for $\la\in X_{*,\R}(T_G)$, one
has $\|\la\|_G=\sqrt{(\la,\la)_{*,G}}$.
This last observation implies that $(.,.)_{*,G}$ is invariant under the action of the Weyl group
$W(T_G):=\n_G(T_G)/T_G$. By polarization, this is equivalent to the fact that $\|.\|_G$ restricted
to $X_{*,\R}(T_G)$ is invariant under $W(T_G)$, and this is obvious from the definition.
\begin{Thm}[Kempf]
\label{kempfi}
Suppose $w\in W$ is a $G$-unstable point. Then, the function $\la\lma \nu_\kappa(\la,w):=\mu_\kappa(\la, w)/\|\la\|_G$ 
on the set of all one
parameter subgroups of $G$ attains a minimal value $m_0\in \Q_{<0}$, and there is a unique parabolic subgroup 
$Q(w)\subset G$,
such that $Q(w)=Q_G(\la)$ for every one parameter subgroup $\la\colon {\Bbb G}_m({\Bbb K})\lra G$ with $\nu(\la,w)=m_0$. 
Moreover, if $\la$ and $\la^\p$ are two indivisible one parameter subgroups with 
$\nu(\la,w)=m_0=\nu(\la^\p,w)$, then there
exists a unique element $u\in {\cal R}_u(Q(w))$, such that 
$\la^\p=u\cdot \la\cdot u^{-1}$.
\end{Thm}
\begin{proof} This is Theorem~2.2 in \cite{Kempf}. It is also proved in \cite{RamRam}, Theorem 1.5.
One has to use (\ref{Sprcomparison}) and (\ref{RRcomparison}) to adapt the formulation to the conventions we use.
Since this theorem plays such a crucial r\^ole in our considerations, we briefly remind the reader of the idea of proof.
Recall Equation (\ref{DeFine}) and the fact that $\|\la\|_G=\|g\cdot\la\cdot g^{-1}\|_G$ for all
$\la\colon {\Bbb G}_m({\Bbb K})\lra G$, $g\in G$. First, for an element $g\in G$,
we search for
\begin{equation}
\label{tominimize}
\min\Bigl\{\,\frac{\mu_\kappa(\la,g\cdot w\cdot g^{-1})}{\|\la\|_G}\,\Bigl|\, \la\in X_*(T_G)\,\Bigr\}.
\end{equation}
Write
$$
{\rm WT}(g\cdot w, T) := \bigl\{\,\chi^g_1,...,\chi^g_{s(g)}\,\bigr\}.
$$
We obtain the linear forms 
\begin{eqnarray*}
l^g_i\colon\q X_{*,{\Bbb \R}}(T_G) &\lra& {\Bbb \R}
\\
\la &\lma& \langle\la, \chi^g_i\rangle_{\Bbb R},\q i=1,...,s(g),
\end{eqnarray*}
on $X_{*,{\Bbb R}}(T_G)$ which are actually defined over $\Q$.
One has now to study the function 
$$
l^g\colon \la\lma \max_{i=1,...,s(g)} l^g_i(\la)
$$
on the norm-one hypersurface $H$ in $X_{*,{\Bbb R}}(T_G)$ where the assumption is that $l$ possesses
a negative value. One then shows that a function like $l^g$ admits indeed a minimum in a unique point $h\in H$.
Moreover, the fact that the $l^g_i$ are defined over $\Q$ grants that the ray ${\Bbb R}_{>0}\cdot h$
contains rational and integral points. See Lemma 1.1 in \cite{RamRam} for this discussion.
Thus, the expression (\ref{tominimize}) agrees with
$l^g(h)$.
\par
Finally, one remarks that $l^g$ depends only on the set of weights ${\rm WT}(g\cdot w, T)$ for which there are only
finitely many possibilities, so that
there is a finite set $\Gamma\subset G$ with
$$
{\rm WT}(g\cdot w, T)\in \Bigl\{\,{\rm WT}(\gamma\cdot w, T)\,|\,\gamma\in \Gamma\,\Bigr\},\q \hbox{for all $g\in G$}.
$$
Thus, we have to show that 
$$
\min_{\gamma\in \Gamma}\min\bigl\{\, l^\gamma(\la)\,|\, \la\in H\,\bigr\}
$$
exists, but this is now clear.
\end{proof}
Let $w$ and $m_0$ be as in the theorem.
We call an indivisible one parameter subgroup $\la\colon{\Bbb G}_m({\Bbb K})\lra G$ with $\nu(\la,w)=m_0$ an \it instability one
parameter subgroup for $w$\rm. Note that, by the theorem, every maximal torus of $Q(w)$ contains a unique instability
one parameter subgroup for $w$.
\begin{Rem}
There is also a canonical parabolic subgroup $Q_{\GL_n({\Bbb K})}(w)$ of $\GL_n({\Bbb K})$
with $Q_{\GL_n({\Bbb K})}(w)\cap G=Q(w)$.
Indeed, if $\la$ is any instability subgroup of $w$, then we set $Q_{\GL_n({\Bbb K})}(w):=Q_{\GL_n({\Bbb K})}(\la)$.
This is well-defined because of the last statement in the theorem.
\end{Rem}
For every maximal torus $T^\p$ of $\GL_n({\Bbb K})$, the given product on $X^*_\R(T)$ induces the pairing 
$(.,.)^*_{T^\p}\colon X^*_\R(T^\p)\times X^*_\R(T^\p)\lra \R$,
$(\chi,\chi^\p)\lma (\chi(g\cdot .\cdot g^{-1}),\chi^\p(g\cdot .\cdot g^{-1}))^*$, where $g\in\GL_n({\Bbb K})$ is an element,
such that $g\cdot T\cdot g^{-1}=T^\p$. Here, the invariance of $(.,.)^*$ under the Weyl group $W(T)$ implies that 
this product does not depend on the choice of $g$. We set $H_G(w):= Q(w)/{\cal R}_u(Q(w))$, and 
$H_{\GL_n({\Bbb K})}(w):=Q_{\GL_n({\Bbb K})}(w)/{\cal R}_u(Q_{GL(W)}(w))$. Now, $\la$ defines an antidominant
character on $H_{\GL_n({\Bbb K})}(w)$ as follows: Let $\ol{T}$ be a maximal torus of $H_{\GL_n({\Bbb K})}(w)$. Under the isomorphism
$L_{\GL_n({\Bbb K})}(\la)\lra H_{\GL_n({\Bbb K})}(w)$ induced by the quotient morphism $\pi\colon Q_{\GL_n({\Bbb K})}(w)\lra H_{\GL_n({\Bbb K})}(w)$, 
there is a unique maximal torus
$T^\p\subset L_{\GL_n({\Bbb K})}(\la)$ mapping onto $\ol{T}$. Then, as we have explained before, there is a scalar
product $(.,.)_{T^\p}^*\colon X^*_\R(T^\p)\times X^*_\R(T^\p)\lra \R$. This provides us with the unique element 
$l_{T^\p}(\la)$, such that $(l_{T^\p}(\la),\chi)^*_{T^\p}=\langle\la,\chi\rangle_\R$ for all
$\chi\in  X^*_\R(T^\p)$. The computation below 
(Example~\ref{CharComp}) shows 
that $l_{T^\p}(\la)$ is indeed a character of $L_{\GL_n({\Bbb K})}(\la)$ and, thus, of $H_{\GL_n({\Bbb K})}(w)$. Call this character
$\chi_0$. Let $T^{\p\p}$ be any other maximal torus of $Q_{\GL_n({\Bbb K})}(w)$. Then, there is an element
$p\in Q_{\GL_n({\Bbb K})}(w)$ with $p\cdot T^{\p}\cdot p^{-1}=T^{\p\p}$. For all one parameter subgroups 
$\widetilde{\la}\colon {\Bbb G}_m({\Bbb K})\lra T^\p$, 
we have
$$
\langle p\cdot \widetilde{\la}\cdot p^{-1}, \chi_0\rangle
\q=\q
\langle \widetilde{\la},\chi_0\rangle
\q=\q
(\widetilde{\la},\la)^*_{T^\p}
\q=\q
(p\cdot \widetilde{\la}\cdot p^{-1},p\cdot\la\cdot p^{-1})^*_{T^{\p\p}},
$$
so that $p\cdot\la\cdot p^{-1}$ and the maximal torus $\ol{T}^\p:=\pi(T^{\p\p})$ yield indeed the same character
$\chi_0$.
\begin{Ex}
\label{CharComp}
Fix integers $0=:n_0<n_1<\cdots<n_s<n_{s+1}:=n$ and $\gamma_1<\cdots<\gamma_{s+1}$
with $\sum_{i=1}^{s+1} \gamma_{i}(n_{i}-n_{i-1})=0$. This defines a one parameter subgroup $\la\colon
{\Bbb G}_m({\Bbb K})\lra \SL_{n}({\Bbb K})$ via
$$
\la(z)\cdot b_j:= z^{\gamma_i}\cdot b_j,\q j=n_{i-1}+1,...,n_i,\ i=1,...,s+1.
$$
Here, $b_1,...,b_n$ is the standard basis for ${\Bbb K}^n$. Then, $L_{\GL_n({\Bbb K})}(\la)\cong \GL_{n_1}({\Bbb K})\times
\GL_{n_2-n_1}({\Bbb K})\times \cdots\times \GL_{n-n_s}({\Bbb K})$, the latter group being embedded as a group of block diagonal
matrices into $\GL_n({\Bbb K})$. One checks that 
$$
l_T(\la)(m_1,...,m_{s+1})\q=\q\det(m_1)^{\gamma_1}\cdot...\cdot \det(m_{s+1})^{\gamma_{s+1}},\q
\forall\ (m_1,...,m_{s+1})\in L_{\GL_n({\Bbb K})}(\la).
$$
\end{Ex}
Let $w\in W\setminus\{0\}$ be an unstable point, and let $Q(w)\subset G$ be the associated parabolic subgroup.
Moreover, choose an instability one parameter subgroup $\la\colon\ {\Bbb G}_m({\Bbb K})\lra G$ for $w$. This yields, in particular,
a flag $W^\bullet\colon\ 0\subsetneq W_1\subsetneq \cdots\subsetneq W_{t}\subsetneq W$. Next, set
$j_0:=\min\{\, j=1,...,t+1\,|\, w\in W_j\,\}$. Then, $w$ defines a point 
$x_\infty\in \P\bigl( (W_{j_0}/W_{j_0-1})^\vee\bigr)$. Let $m_0\in \Q_{<0}$ be as in Theorem~\ref{kempfi}, 
and $q:=m_0\cdot \|\la\|_G=\mu_\kappa(\la,w)\in\Z_{<0}$. Finally, define $\chi_*:= q\cdot \chi_{0|H_G(w)}$.
\begin{Prop}[Ramanan-Ramanathan]
\label{RamRamI}
The point $x_\infty\in \P\bigl( (W_{j_0}/W_{j_0-1})^\vee\bigr)$ is semistable for the induced $H_G(w)$-action
and its linearization in $\O_{\P( (W_{j_0}/W_{j_0-1})^\vee)}(1)$ twisted by the character $\chi_*$.
\end{Prop}
\begin{proof} This is Proposition 1.12 in \cite{RamRam}. We observe that, by 
(\ref{Sprcomparison}) and (\ref{RRcomparison}),
we have $\chi_*=\chi$ with $\chi$ the character constructed in  \cite{RamRam}.
(Our explicit construction shows that we may take $s=1$ and $r=1$ in the proof of \cite{RamRam}, Proposition 1.12).
Note that Ramanan and Ramanathan
show that $x_\infty=[w_\infty]$, where 
$w_\infty\in W_{j_0}/W_{j_0-1}\otimes {\Bbb K}_{\chi_*^{-1}}$ is a semistable point and ${\Bbb K}_{\chi_*^{-1}}$ 
is the
one dimensional $H_G(w)$-module associated with the character $\chi_*^{-1}$. This gives the claimed linearization.
\end{proof}
Finally, we need Kempf's rationality result. For this, let $K$ be a non-algebrai\-cally closed field of characteristic
zero (in our application, this will be the function field of an algebraic variety), $G\lra \Spec(K)$ a $K$-group,
and $W$ a finite dimensional $K$-vector space. Fix an algebraic closure ${\Bbb K}$ of $K$, and set $G_{\Bbb K}
:= G\times_{\Spec(K)}\Spec(\Bbb K)$ and $W_{\Bbb K}:=W\otimes_K{\Bbb K}$. Suppose that we are
given a $K$-rational representation $\kappa\colon\ G_{\Bbb K}\lra \GL(W_{\Bbb K})$.
\begin{Thm}[Kempf]
\label{rationality}
If $T\subset G_{\Bbb K}$ is a maximal torus which is defined over $K$ and $K_T/K$ is a finite extension
of $K$, such that $T_K\times_{\Spec(K)}\Spec(K_T)\cong {\Bbb G}_a(K_T)\times\cdots\times{\Bbb G}_a(K_T)$,
$T_K\subset G$ being the $K$-group with $T_K\times_{\Spec(K)}\Spec(\Bbb K)=T$, then, for a product
$(.,.)^*\colon X_\R^*(T)\times X^*_\R(T)\lra \R$ which is invariant under both the Weyl group $W(T)$ and
the action of the Galois group ${\rm Gal}({\Bbb K}/K)$ via its finite quotient ${\rm Gal}(K_T/K)$, 
the following holds true: 
If $w\in W_{\Bbb K}$ is an unstable $K$-rational point, then the
parabolic subgroup $Q_{G_{\Bbb K}}(w)$, associated to $w$ by means of the norm $\|.\|_{G_{\Bbb K}}$ 
on the one parameter subgroups of $G_{\Bbb K}$ which is induced by $(.,.)^*$, is defined over $K$.
\end{Thm}
\begin{proof}
This is part of Theorem 4.2 in \cite{Kempf}.
See also \cite{RamRam} for generalizations.
\end{proof}
\paragraph{Weighted Projective Spaces. ---}
For a given tuple $\ul{w}=(w_0,...,w_n)$ of positive integers, the quotient of 
${\Bbb A}^{n+1}\setminus\{0\}$ w.r.t.\ to the $\C^*$-action $z\cdot (x_0,....,x_n)=
\bigl(z^{w_0}\cdot x_0,...,z^{w_n}\cdot x_n\bigr)$ is the so-called \it weighted projective space \rm
$\P_{\ul{w}}$. One has $\P_{\ul{w}}\cong {\rm Proj}(\C[y_0,...,y_n])$ where one assigns the weight 
$w_i$ to the variable
$y_i$, $i=0,...,n$. Then, the degree is defined for every monomial in the $y_i$, and for each non-negative integer
$\omega$, we can speak of the homogeneous elements of degree $\omega$ and define
$\C[y_0,...,y_n]^{\omega}\subset \C[y_0,...,y_n]$ as the finite dimensional vector space generated by the homogeneous
elements of degree $\omega$. We also define the subalgebra
$$
\C[y_0,...,y_n]^{(\omega)}:=\bigoplus_{i\in \Z_{\ge 0}} \C[y_0,...,y_n]^{i\omega}.
$$
If $\omega$ is a sufficiently large common multiple of $w_0,...,w_n$, the subalgebra
$\C[y_0,...,y_n]^{(\omega)}$ is generated by ${\Bbb V}_\omega:=\C[y_0,...,y_n]^{\omega}$ (\cite{RedBook},
Chapter III, p.\ 282), i.e., we have
a surjection $\Sym^*{\Bbb V}_\omega\lra \C[y_0,...,y_n]^{(\omega)}$ that defines an embedding 
$\iota\colon \P_{\ul{w}}\hookrightarrow \P({\Bbb V}_\omega)$.
\begin{Rem}
Alternatively, pick primitive
$w_i$-th roots of unity $\zeta_i$, $i=1,...,n$, and look at the action of $\Z_{\ul{w}}:=\Z_{w_0}\oplus\cdots
\oplus \Z_{w_n}$ on $\P_n$ by 
$$
(b_0,...,b_n)\cdot [x_0:\cdots:x_n]\ :=\ \bigl[\zeta_0^{b_0}\cdot x_0:\cdots:\zeta_n^{b_n}\cdot x_n\bigr],
\q (b_0,...,b_n)\in \Z_{\ul{w}}, [x_0:\cdots:x_n]\in\P_n.
$$
Then, $\P_n/\Z_{\ul{w}}\cong\P_{\ul{w}}$, where the quotient morphism 
$\P_n\lra \P_{\ul{w}}$ corresponds to the algebra
homomorphism $\C[y_0,...,y_n]\lra \C[x_0,...,x_n]$, $y_i\lma x_i^{w_i}$, $i=0,...,n$. 
\end{Rem}
We will need an intrinsic description of the morphism 
$\psi\colon {\Bbb A}^{n+1}\setminus\{0\}\lra \P({\Bbb V}_\omega)$. 
For this,
let $v_1<\cdots<v_m$ be the different integers appearing in $\ul{w}$, and let
$$
\C^{n+1}\ =\  V_1\oplus\cdots\oplus V_m
$$
be the corresponding decomposition into eigenspaces. For $\omega$ as above, the morphism $\psi$ comes
from
$$
\bigoplus_{{(d_1,...,d_m):}\atop{v_1d_1+\cdots+v_md_m=\omega}} S^{d_1}\bigl(V_1^\vee\bigr)
\otimes \cdots\otimes S^{d_m}\bigl(V_m^\vee\bigr)\otimes\O_{{\Bbb A}^{n+1}\setminus\{0\}}
\stackrel{\sum S^{d_1+\cdots+d_m}(\tau)}
{\lra}
\O_{{\Bbb A}^{n+1}\setminus\{0\}}
$$
with $\tau\colon {\C^{n+1}}^\vee\otimes\O_{{\Bbb A}^{n+1}\setminus\{0\}}\lra \O_{{\Bbb A}^{n+1}\setminus\{0\}}$ 
the dual of the tautological section.
\par
One checks easily that
\begin{equation}
\label{equality501}
\mu_{\id}(\la, v)\ <\ (=/>)\ 0
\q\Longleftrightarrow\q \mu_{\O_{\P({\Bbb V}_\omega)}(1)}\bigl(\la, \psi(v)\bigr)\ <\ (=/>)\ 0,
\end{equation}
for every point $v\in{\Bbb A}^{n+1}\setminus\{0\}$ and every one parameter subgroup 
$\la\colon \C^*\lra\GL_{n+1}(\C)$.
\subsection{Principal Bundles}
Let $U$ be a smooth algebraic variety and $G$ a reductive algebraic group over the field of complex numbers. 
Suppose we are given a principal
$G$-bundle ${\cal P}$ over $U$. If $F$ is an algebraic variety and $\alpha\colon G\times F\lra F$ is an 
action of $G$ on $F$, then we may form the geometric quotient
$$
{\cal P}(F,\alpha):= \bigl({\cal P}\times F\bigr)/G
$$ 
w.r.t.\ the action $(p,f)\cdot g:= (p\cdot g, g^{-1}\cdot f)$ for all $p\in {\cal P}$, $f\in F$, and $g\in G$. Note that
${\cal P}(F,\alpha)$ is a fibre space with fibre $F$ over $U$ which is
locally trivial in the \'etale topology. An important 
special case arises when we look at the action $\frak c\colon G\times G\lra G$, $(g,h)\lma g\cdot h\cdot g^{-1}$,
of $G$ on itself by conjugation. Then, the associated fibre space ${\cal G}({\cal P}):={\cal P}(G,\frak c)\lra U$ is 
a reductive group scheme over $U$, and, for any pair $(F,\alpha)$ as above, we obtain an induced action
$$
a\colon\ {\cal G}({\cal P})\times_U {\cal P}(F,\alpha)\lra {\cal P}(F,\alpha).
$$
If $W$ is a vector space and $\kappa\colon G\lra \GL(W)$ is a representation, we set ${\cal P}_\kappa:={\cal P}(W,\kappa)$.
Note that the formation of ${\cal P}(F,\alpha)$ commutes with base change. For additional information, we refer
the reader to \cite{Serre}.
\paragraph{Parabolic Subgroup Schemes. ---}
Let $S$ be any scheme and
suppose ${\cal G}_S\lra S$ is a reductive group scheme over $S$. A subgroup ${\cal Q}_S\subset {\cal G}_S$ is called
a \it parabolic subgroup\rm, if it is smooth over $S$ and, for any geometric point $s$ of $S$, the quotient
${\cal G}_{S,s}/{\cal Q}_{S,s}$ is proper. The functor
\begin{eqnarray*}
\ul{\rm Par}({\cal G}_S)\colon\ \ul{\rm Schemes}_S &\lra& \ul{\rm Sets}
\\
(T\lra S) &\lma& \Bigl\{\, \hbox{Parabolic subgroups of ${\cal G}_S\times_S T$}\,\Bigr\}
\end{eqnarray*}
is then representable by an $S$-scheme ${\cal P}ar({\cal G}_S)$. For the details, we refer the reader to \cite{Dem}.
\begin{Ex}
\label{ParRep}
Let $G$ be a complex reductive group and $\cal P$ be a principal $G$-bundle over the variety $U$. 
Denote by $\frak P$ the
set of conjugacy classes of parabolic subgroups in $G$ and pick for each class $\frak p\in \frak P$ a representative
$Q_\frak p$. Then,
$$
{\cal P}ar\bigl({\cal G}({\cal P})\bigr)\cong \bigsqcup_{\frak p\in\frak P} {\cal P}/Q_\frak p.
$$
\end{Ex}
\paragraph{Sections in Associated Projective Bundles. ---}
Let ${\cal P}\lra U$ be a principal bundle as before. Suppose we are given a representation $\kappa\colon
G\lra \GL(W)$. This yields an action $\alpha\colon G\times \P(W^\vee)\lra \P(W^\vee)$ and a linearization
$\ol{\alpha}\colon G\times \O_{\P(W^\vee)}(1)\lra \O_{\P(W^\vee)}(1)$ of this action. 
Let $\P^{\rm ss}\subset \P(W^\vee)$ be the open subset of semistable points.
\begin{Prop}[Ramanan-Ramanathan] 
\label{RamRamII}
Assume that $U$ is a big open subset of the manifold $X$. Let  
$\sigma\colon U\lra {\cal P}(\P(W^\vee),\alpha)$ be a section, and $\L_\sigma$ the pullback --- via
$\sigma$ --- of
the line bundle ${\cal P}(\O_{\P(W^\vee)}(1),\ol{\alpha})\lra {\cal P}(\P(W^\vee),\alpha)$ to $U$. 
If $\sigma(\eta)\in
{\cal P}(\P^{\rm ss},\alpha)$, then
$$
\deg(\L_\sigma)\q\ge\q 0.
$$
Here, $\eta$ is the generic point of $U$.
\end{Prop}
\begin{proof} This is Proposition 3.10, i), in \cite{RamRam}.
\end{proof}
\section{Inhomogeneous Decorations}
\label{DecSh}
\setcounter{subsection}{1}
Fix tuples $\ul{a}=(a_1,...,a_n)$, $\ul{b}=(b_1,...,b_n)$, and $\ul{c}=(c_1,...,c_n)$ of non-negative integers,
such that $a_i-rc_i>0$, for $i=1,...,n$. 
If $a_i-rc_i=a_j-rc_j$, $1\le i<j\le n$, we call the triple $(\ul{a},\ul{b},\ul{c})$ \it homogeneous\rm.
We look at pairs $(\A,\phi)$
with $\A$ a torsion free sheaf of rank $r$ and $\phi\colon \A_{\ul{a},\ul{b},\ul{c}}\lra\O_X$ a non-trivial
homomorphism. Here, 
$$
\A_{\ul{a},\ul{b},\ul{c}}\q:=\q \bigoplus_{i=1}^n 
\bigl(\A^{\otimes a_i}\bigr)^{\oplus b_i}\otimes \bigl(\det(\A)^\vee\bigr)^{\otimes c_i}.
$$
We call $(\A,\phi)$ a \it torsion free sheaf with a decoration of type $(\ul{a},\ul{b},\ul{c})$\rm, and we
say that $(\A,\phi)$ is \it equivalent to \rm $(\A^\p,\phi^\p)$, if there is an isomorphism
$\psi\colon \A\lra\A^\p$, such that $\phi=\phi^\p\circ \psi_{\ul{a},\ul{b},\ul{c}}$, letting
$\psi_{\ul{a},\ul{b},\ul{c}}\colon \A_{\ul{a},\ul{b},\ul{c}}\lra \A^\p_{\ul{a},\ul{b},\ul{c}}$
be the isomorphism induced by $\psi$. The decoration
$\phi$ breaks into components
$$
\phi_i\colon\ \bigl(\A^{\otimes a_i}\bigr)^{\oplus b_i}\lra \det(\A)^{\otimes c_i},\q i=1,...,n.
$$
Given a weighted filtration $(\A^\bullet, \ul{\alpha})$, we define the weight vector 
$$
\ul{\gamma}=(\gamma_1,...,\gamma_{s+1})\q:=\q
\sum_{i=1}^s\alpha_i\cdot\bigl(\underbrace{\rk(\A_i)-r,...,\rk(\A_i)-r}_{\rk(\A_i)\times},
\underbrace{\rk(\A_i),...,\rk(\A_i)}_{(r-\rk(\A_i))\times}\bigr),\  r:=\rk(\A),
$$ 
and set, for $\A_{s+1}:=\A$ and $i=1,...,n$,
$$
\mu\bigl(\A^\bullet,\ul{\alpha};\phi_i\bigr):= 
-\min\Bigl\{\gamma_{j_1}+\cdots+\gamma_{j_{a_i}}\,|\, {(j_1,...,j_{a_i})\in\{\,1,...,s+1\,\}^{\times a_i}}
:\phi_{i|(\A_{j_1}\otimes\cdots\otimes\A_{j_{a_i}})^{\oplus b_i}}\not\equiv 0\Bigr\}
$$
as well as
$$
\mu\bigl(\A^\bullet,\ul{\alpha};\phi\bigr)\q:=\q\max\Bigl\{\,
\mu\bigl(\A^\bullet,\ul{\alpha};\phi_i\bigr)\,|\, i=1,...,n\,\Bigr\}.
$$
\begin{Ex}
\label{weightcomputation}
Let $G$ be a reductive group, $\rho\colon G\lra\GL(V)$ a faithful representation, and 
$\kappa\colon \GL(V)\lra \GL(W)$ a representation of $\GL(V)$.
For non-negative integers $a,b,c$, we have the $\GL(V)$-module
$$
V_{a,b,c}\q=\q \bigl(V^{\otimes a}\bigr)^{\oplus b}\otimes \bigl(\bigwedge^{r} V\bigr)^{\otimes -c},
\q r:=\dim(V),
$$
and, for tuples $\ul{a}=(a_1,...,a_n)$, $\ul{b}=(b_1,...,b_n)$, and $\ul{c}=(c_1,...,c_n)$, 
we define the $\GL(V)$-module
$$
V_{\ul{a},\ul{b},\ul{c}}\q:=\q \bigoplus_{i=1}^n V_{a_i,b_i,c_i}.
$$
The corresponding representation is denoted by $\kappa_{\ul{a},\ul{b},\ul{c}}$. We may assume
that $W$ is a submodule of $V_{\ul{a},\ul{b},\ul{c}}$ (\cite{Schmitt0}, Section 1.1).
Thus, there exists a $G$-submodule $W^\p$ of $V_{\ul{a},\ul{b},\ul{c}}$ which is complementary to
$W$ and we find an isomorphism $V_{\ul{a},\ul{b},\ul{c}}\cong W\oplus W^\p$ of $G$-modules. 
Furthermore, suppose $(\A,\tau)$ is an honest singular principal $G$-bundle. Then, over the open set
$U$ where $\A$ is locally free, we have the principal bundle ${\cal P}:={\cal P}(\A,\tau)$.
Viewing $\kappa$ as a representation of $G$, we have the associated vector bundle 
${\cal P}_{\kappa}$ with fibre $W$. By our assumption,
this is a subbundle of ${\A}^\vee_{|U;\ul{a},\ul{b},\ul{c}}$. Thus,  a section
$$
\phi^\p\colon\ U\lra {\cal P}_{\kappa}
$$
defines components 
$$
\phi^\p_i\colon\ U\lra \Bigl({\A^\vee_{|U}}^{\otimes a_i}\Bigr)^{\oplus b_i}\otimes \det(\A_{|U})^{\otimes c_i},\q
i=1,...,n,
$$
which correspond to
$$
\phi^{\p\p}_i\colon \bigl(\A_{|U}^{\otimes a_i}\bigr)^{\oplus b_i}\lra \det(\A_{|U})^{\otimes c_i},\q
i=1,...,n,
$$
and the latter extend uniquely to
$$
\phi_i\colon \bigl(\A^{\otimes a_i}\bigr)^{\oplus b_i}\lra \det(\A)^{\otimes c_i},\q
i=1,...,n.
$$
Altogether, we obtain an associated torsion free sheaf $\A$ with a decoration $\phi$ of type $(\ul{a},\ul{b},\ul{c})$.
If we are given a one parameter subgroup $\la\colon\C^*\lra G$ and a reduction $\beta$ of $(\A,\tau)$
to the one parameter subgroup $\la\colon\C^*\lra G$, then we can give
a more intuitive description of  $\mu(\A_\beta^\bullet,\ul{\alpha}_\beta;\phi)$. Over the open subset $U^\p$, 
the reduction $\beta$ gives rise to a $Q_G(\la)$-bundle.
Making $U$ smaller and passing to an \'etale covering $\widetilde{U}$, 
we may trivialize
this $Q_G(\la)$-bundle. Such a trivialization gives rise to an isomorphism
$$
\Psi_{\widetilde{U}}\colon\ \A^\vee_{|\widetilde{U}}\ \cong\  V\times \widetilde{U}
$$
with
$$
\Psi_{\widetilde{U}}(\A^{\p}_i)\ =\  V_i\times \widetilde{U},\q i=1,...,s.
$$
Here, $( 0\subsetneq V_1\subsetneq \cdots\subsetneq V_s\subsetneq
V,\ul{\alpha}\bigr)$ is the weighted flag of $\la$ and 
the $\A^{\p}_i$, $i=1,...,s$, are as in the introduction. Looking
at
$$
\phi^\p_{|\widetilde{U}}\colon\ \widetilde{U}\lra {\cal P}_{\kappa|\widetilde{U}}\cong
W\times\widetilde{U}\lra W,
$$
one finally has 
\begin{equation}
\label{weightcomp2}
\max\Bigl\{\,\mu_\kappa\bigl(\la,
\phi^\p_{|\widetilde{U}}(x)\bigr)\,|\, x\in \widetilde{U}\,\Bigr\}=
\mu(\A_\beta^\bullet,\ul{\alpha}_\beta;\phi).
\end{equation}
This is explained in detail in Section 2.1.1 of \cite{Schmitt0}.
\end{Ex}
Let $v_1<\cdots<v_m$ be the integers which occur as $a_i-rc_i$, $i=1,...,n$. Set $V_j:=\bigoplus_{i:a_i-rc_i=v_j}
V_{a_i,b_i,c_i}$, $j=1,...,m$. Choose a sufficiently large common multiple $\omega$ of $v_1,...,v_m$ as in the section 
on weighted projective spaces. Then, letting $\C^*$ act on $V_j$ by multiplication with $z^{v_j}$,
the weighted projective space $(V_{\ul{a},\ul{b},\ul{c}}\setminus\{0\})/\C^*$ gets embedded into 
$\P({\Bbb V}_\omega)$,
$$
{\Bbb V}_\omega\q:=\q\bigoplus_{(d_1,...,d_m):\atop v_1d_1+\cdots+v_md_m=\omega} 
S^{d_1}(V^\vee_1)\otimes\cdots\otimes
S^{d_m}(V^\vee_m).
$$
We may find positive integers $A,B,C$ with $A-rC=-\omega$, such that ${\Bbb V}_\omega$ is a direct summand of $V_{A,B,C}$
(\cite{Schmitt0}, Corrolary 1.2)
and we have an embedding $\P({\Bbb V}_\omega) \hookrightarrow\P(V_{A,B,C})$.
\par
Let $({\A},\phi)$ be a torsion free sheaf with a decoration of type $(\ul{a},\ul{b},\ul{c})$, then $\phi$ defines,
for every tuple $\ul{d}=(d_1,...,d_m)$ with $v_1d_1+\cdots+v_md_m=\omega$, a homomorphism 
$$
\widetilde{\phi}_{\ul{d}}\colon\  S^{d_1}({\cal B}_1)\otimes\cdots\otimes
S^{d_m}({\cal B}_m)\lra \O_{X}.
$$
These homomorphisms add to
$$
\widetilde{\phi}\colon\ \bigoplus_{(d_1,...,d_m):\atop v_1d_1+\cdots+v_md_m=\omega} 
S^{d_1}({\cal B}_1)\otimes\cdots\otimes
S^{d_m}({\cal B}_m)\lra \O_X,
$$
${\cal B}_j:=\bigoplus_{i:a_i-rc_i=v_j} \A_{a_i,b_i,c_i}$, $j=1,...,m$. This finally defines
$$
\widehat{\phi}\colon\ \bigl(\A^{\otimes A}\bigr)^{\oplus B}\lra \det(\A)^{\otimes C}.
$$
For every weighted filtration $(\A^\bullet,\ul{\alpha})$ of $\A$, we set
$$
\nu(\A^\bullet,\ul{\alpha};\phi)\q:=\q\frac{1}{\omega}\cdot \mu(\A^\bullet,\ul{\alpha};\widehat{\phi}).
$$
Next, we also fix a positive polynomial $\eps\in\Q[x]$ of degree at most $\dim(X)-1$. Then, a torsion free
sheaf $(\A,\phi)$ with a decoration of type $(\ul{a},\ul{b},\ul{c})$ is called \it $\eps$-(semi)stable\rm, if,
for every weighted
filtration $(\A^\bullet,\ul{\alpha})$, the inequality
$$
M(\A^\bullet,\ul{\alpha})+\eps\cdot \nu(\A^\bullet,\ul{\alpha};\phi)\q(\succeq)\q 0
$$
is satisfied. This gives rise to a moduli problem which is reduced to the homogeneous case ---
covered by \cite{Schmitt0} and \cite{GS} --- via the above assignment $\phi\lma\widehat{\phi}$
and then solved along the lines of our construction in \cite{Schmitt}
and will, therefore, not be explained here. The crucial fact is that, for a given Hilbert polynomial $P$
and a given polynomial $\eps$ as above, the $\eps$-semistable pairs $(\A,\phi)$ with $P({\cal A})=P$
form a \sl projective \rm moduli space.
\par
An important fact is that (\ref{equality501}) implies
\begin{equation}
\label{equality502}
\nu(\A^\bullet,\ul{\alpha};\phi)\ <\ (=/>)\ 0
\q\Longleftrightarrow\q \mu(\A^\bullet,\ul{\alpha};\phi)\ <\ (=/>)\ 0.
\end{equation}
\begin{Rem}
\label{MarB}
i) Unlike the situation in \cite{Schmitt}, we can choose $\omega$ here to work for all the objects we consider,
because it obviously depends only on the fixed input data $\ul{a}$, $\ul{b}$, and $\ul{c}$.
\par
ii) Let $(\A,\widehat{\phi})$ be as above. Then, as in \cite{Schmitt0}, Lemma 1.8, i), 
for every saturated subsheaf $0\subsetneq {\cal B}\subsetneq \A$, one finds
$\mu(0\subsetneq {\cal B}\subsetneq \A, (1); \widehat{\phi})\le A(r-1)$. Thus, the semistability condition
yields
$$
\mu({\cal B}) \le \mu(\A)+\frac{\eps_0\cdot A\cdot (r-1)}{r\cdot \omega},
$$
$\eps_0$ being the coefficient of the monomial of degree $\dim X-1$ in $\eps$.
\end{Rem}
The following result will be the basis of our ``semistable reduction theorem'' for semista\-ble honest principal
$G$-bundles, and, thus, constitutes a core result of this paper. Its proof follows the strategy of 
Ramanan and Ramanathan
in their proof of Proposition~3.13.
\begin{Thm}
\label{SemStabRed}
Suppose the stability parameter $\eps$ has degree exactly $\dim(X)-1$ and that $(\A,\phi)$ is an $\eps$-semistable
torsion free sheaf with a decoration of type $(\ul{a},\ul{b},\ul{c})$, such that $\deg(\A)=0$. 
Then the following holds true: Denote by $\eta$ the generic point of $X$ and by $K$ its residue field, 
and choose a trivialization
$\A^\vee\otimes_{\O_X} \O_{X,\eta}\cong V\times_{\Spec(\C)} \Spec(K)$.
Then, the point $\widetilde{\sigma}_\eta\in V_{\ul{a},\ul{b},\ul{c}}\times_{\Spec(\C)} \Spec(K)$ defined by $\phi$ 
and the trivialization lies
in $V^{\rm ss}_{\ul{a},\ul{b},\ul{c}}\times_{\Spec(\C)} \Spec(K)$. Here, $V^{\rm ss}_{\ul{a},\ul{b},\ul{c}}$
is the open subset of $\SL(V)$-semistable points in $V_{\ul{a},\ul{b},\ul{c}}$.
\end{Thm}
\begin{proof}
Let us start with some notation. We use the representation $\kappa:=\kappa_{\ul{a},\ul{b},\ul{c}}
\colon \GL(V)\lra \GL(V_{\ul{a},\ul{b},\ul{c}})$ and let 
$\alpha\colon \GL(V)\times Y\lra Y$, $Y:=\P(V^\vee_{\ul{a},\ul{b},\ul{c}})$, be the
induced action. Moreover, ${\cal P}:=\ul{\rm Isom}(V\otimes\O_U,\A_{|U}^\vee)$ is the principal $\GL(V)$-bundle
associated with $\A_{|U}$ over the maximal open subset $U$ over which $\A$ is locally free. 
The group $\GL(V)$ acts on $\SL(V)$ and $\GL(V)$ by conjugation $\frak c$, and we let
${\cal SL}({\cal P})\subset {\cal GL}({\cal P})$ be the corresponding group schemes over $U$.
We fix an algebraic
closure ${\Bbb K}$ of $K$. A trivialization as chosen in the statement of the theorem
is equivalent to a trivialization ${\cal P}\times_U \Spec(K)\cong \GL(V)\times_{\Spec(\C)} \Spec(K)$.
The latter identification will induce trivializations of the objects introduced below.
Define
$$
\ol{Y}:= {\cal P}(Y,\alpha)\times_U \Spec({\Bbb K})\cong Y\times_{\Spec(\C)} \Spec({\Bbb K}),
$$
$G:=\SL(V)$, and
$$
\ol{G}:= {\cal SL}({\cal P})\times_U \Spec({\Bbb K})\cong G\times_{\Spec(\C)} \Spec({\Bbb K}).
$$
Finally, set $L:=\O_Y(1)$, and $\ol{L}:=\O_{\ol{Y}}(1)\cong L\times_{\Spec(\C)} \Spec({\Bbb K})$.
\par
Next, we remind the reader of Proposition 1.14 in Chapter 1.4 of \cite{GIT}:
\begin{Prop}[Mumford]
\label{Mummy}
The set of $\ol{G}$-semistable points in $\ol{Y}$ w.r.t.\ the linearization in $\ol{L}$ is given as
$$
Y^{\rm ss}(L)\times_{\Spec(\C)} \Spec({\Bbb K}).
$$
Here, $Y^{\rm ss}(L)$ is the set of $G$-semistable points in $Y$ w.r.t.\ linearization in $L$.
\end{Prop}
Now, let ${\cal NC}\subset V_{\ul{a},\ul{b},\ul{c}}$ be the cone of $\SL(V)$-unstable points.
Recall that we are given a section $\sigma_U\colon U\lra {\cal P}_\kappa$. Set
$\sigma_\eta:=\sigma_U\times_U \Spec({\Bbb K})\in \ol{Y}$. The negation of the assertion of the theorem is,
by Proposition~\ref{Mummy},
$$
\sigma_\eta\ \in\ {\cal NC}\times_{\Spec(\C)} \Spec({\Bbb K}).
$$
Our first step toward the proof will be an application of Kempf's rationality theorem \ref{rationality}.
For this, let $T\subset\GL(V)$ be a maximal torus. We may choose a basis of $V$, such that $T$ becomes
the subgroup of diagonal matrices. Then, we define the pairing $(.,.)^*\colon X^*_{\Bbb R}(T)\times X^*_{\Bbb R}(T)
\lra \R$ as in Section \ref{GIT} (The Instability Flag). Now, $\ol{T}:= T\times_{\Spec(\C)} \Spec({\Bbb K})$
is a maximal torus in $\GL(V)\times_{\Spec(\C)} \Spec({\Bbb K})$ with $X^*(\ol{T})=X^*(T)$, 
and its intersection $T_{\ol{G}}$ with
$\ol{G}$ is a maximal torus in that group. The induced pairing $(.,.)^*_G$ on $X_\R^*(T_{\ol{G}})$ fulfills
the requirements of Theorem~\ref{rationality}. If we assume that $\sigma_\eta$ be unstable, then
there is an instability one parameter subgroup $\la\colon {\Bbb G}_m({\Bbb K})\lra \ol{G}$ which
defines a weighted flag $(0\subsetneq \ol{V}_1\subsetneq\cdots\subsetneq \ol{V}_s\subsetneq \ol{V},\ul{\alpha}^\p)$
in $\ol{V}:=V\otimes_\C {\Bbb K}$. The resulting parabolic subgroup $Q_{\ol{G}}(\sigma_\eta)$ is 
defined over $K$, i.e., it comes from a parabolic subgroup $Q_{K}(\sigma_\eta)$ of ${\cal SL}({\cal P})\times_U
\Spec(K)$. The parabolic subgroup $Q_{K}(\sigma_\eta)$, in turn, corresponds to a point
$$
\Spec(K)\lra {\cal P}ar({\cal SL}({\cal P}))={\cal P}ar({\cal GL}({\cal P}))\stackrel{{\rm \ref{ParRep}}}{\cong}
{\cal P}/Q_{\frak p},
$$
for the appropriate conjugacy class ${\frak p}\in {\frak P}$ of parabolic subgroups of $\GL(V)$.
This point may be extended to a section
$$
U^\p\lra 
{\cal P}/Q_{\frak p}
$$
over a non-empty open subset $U^\p\subset U$. In fact, we may assume $U^\p$ to be big. This is because
$X$ and ${\cal P}/Q_{\frak p}$ are smooth projective varieties, so that any rational map
$X\dasharrow {\cal P}/Q_{\frak p}$ extends to a big open subset. As in the introduction,
this defines a filtration 
$$
0\subsetneq {\cal A}^\p_1\subsetneq \cdots \subsetneq {\cal A}^\p_s\subsetneq {\cal A}^\vee_{|U^\p}
$$
of ${\cal A}^\vee_{|U^\p}$ by subbundles,
and leads to a filtration 
$(\A^\bullet\colon 0\subsetneq \A_1\subsetneq \cdots \subsetneq \A_s\subsetneq \A)$
of $\A$ by saturated subsheaves. We define the vector $\ul{\alpha}=(\alpha_1,...,\alpha_s)$
by $\alpha_i:=\alpha^\p_{s+1-i}$, $i=1,...,s$, if $\ul{\alpha}^\p=(\alpha^\p_1,...,\alpha_s^\p)$.
For the weighted filtration $(\A^\bullet, \ul{\alpha})$, we clearly find
$$
\mu\bigl(\A^\bullet, \ul{\alpha}; \phi\bigr)\q <\q 0.
$$
Thus, by (\ref{equality502}), we infer that $\eps\cdot \nu(\A^\bullet, \ul{\alpha}; \phi)$ is a negative
polynomial of degree exactly $\dim(X)-1$. 
The following claim settles the theorem.
\begin{Claim}
For the weighted filtration $(\A^\bullet, \ul{\alpha})$ constructed above, the coefficient of the monomial 
of degree ${\dim(X)-1}$
in $M(\A^\bullet,\ul{\alpha})$ is not positive.
\end{Claim}
We now explain the proof of this claim.
Note that we obtain, in fact, an even stronger rationality theorem. The group $Q_{\frak p}$ is the stabilizer
of a unique flag $V^\bullet\colon 0\subsetneq V_1\subsetneq\cdots \subsetneq V_s\subsetneq V$, and
the weighted filtration $(\A^\bullet, \ul{\alpha})$ defines a reduction of the structure group of ${\cal P}$ to 
$Q_{\frak p}$. If we start our arguments with a trivialization of the induced 
$Q_{\frak p}$-bundle ${\cal Q}$
over the generic point $\eta$, then we get $\ol{V}^\bullet$ with $\ol{V}_i=V_i\otimes_\C {\Bbb K}$, $i=1,...,s$,
as the instability flag. One may use the weighted flag 
$(V^\bullet, \ul{\alpha}^\p)$ to define
a one parameter subgroup $\la\colon \C^*\lra G$ (which, indeed, is an instability subgroup). 
Then, $\la$ defines also a flag 
$W^\bullet\colon\ 0\subsetneq W_1\subsetneq \cdots\subsetneq W_{t}
\subsetneq W$ in $W:=V_{\ul{a},\ul{b},\ul{c}}$, 
and the parabolic subgroup $Q_{G}(\la)\subset Q_{\GL(V)}(\la)=Q_{\frak p}$ fixes this flag. 
Recall that we are given a reduction of the structure group of ${\cal P}_{|U^\p}$ 
to $Q_{\GL(V)}(\la)=Q_{\frak p}$. Therefore, the flag $W^\bullet$ gives rise to a filtration 
$$
0\subsetneq {\cal B}^\vee_1\subsetneq \cdots\subsetneq {\cal B}^\vee_{t}\subsetneq \A_{|U^\p;\ul{a},\ul{b},\ul{c}}^\vee
$$
by subbundles.
Define 
$$
j_0:=\min\{\, j=1,...,t+1\,|\, {\cal B}^\vee_j\hbox{ contains the image of $\sigma_{U}$}\,\}.
$$
Let $\L^\p\subset \O_X$ be the image of $\phi$. Then, over a big open subset $U^{\p\p}\subset U^\p$,
we have $\L^\p_{|U^{\p\p}}\cong \O_{U^{\p\p}}(-D)$ for an effective divisor $D$. Thus,
$\phi_{|U^{\p\p}}\colon \A_{|U^{\p\p}; \ul{a},\ul{b},\ul{c}}\lra \O_{U^{\p\p}}(-D)$ defines a morphism
$$
\ol{\sigma}_{U^{\p\p}}\colon U^{\p\p}\lra \P(\A_{|U^{\p\p}; \ul{a},\ul{b},\ul{c}})
$$
with
$$
{\ol{\sigma}^*_{U^{\p\p}}}\bigl(\O_{\P(\A_{|U^{\p\p}; \ul{a},\ul{b},\ul{c}})}(1)\bigr)\cong\O_{U^{\p\p}}(-D).
$$
By our choice of $j_0$, 
$\ol{\sigma}_{U^{\p\p}}$ factorizes over $\P({\cal B}_{j_0|U^{\p\p}})$,
and, again,
$$
{\ol{\sigma}^*_{U^{\p\p}}}\bigl(\O_{\P({\cal B}_{j_0|U^{\p\p}})}(1)\bigr)\cong\O_{U^{\p\p}}(-D).
$$
Now, the  surjective linear map $W_{j_0}\lra W_{j_0}/W_{j_0-1}$ is, in fact, a morphism
of $Q_{\frak p}$-modules.
Over a big open subset $U^{\p\p\p}\subset U^{\p\p}$, the image of 
$$
\bigl({\cal B}^\vee_{j_0|U^{\p\p\p}}/{\cal B}^\vee_{j_0-1|U^{\p\p\p}}\bigr)^\vee \subset {\cal B}_{j_0|U^{\p\p\p}}
\lra \O_{U^{\p\p\p}}(-D)
$$
is of the form $\O_{U^{\p\p\p}}(-(D_{|U^{\p\p\p}}+D^\p))$ for some effective divisor $D^\p$. Therefore, we get
a morphism
$$
\ol{\sigma}^{\p\p}\colon U^{\p\p\p}\lra 
\P\bigl(({\cal B}^\vee_{j_0|U^{\p\p\p}}/{\cal B}^\vee_{j_0-1|U^{\p\p\p}})^\vee\bigr)
$$
with
$$
{\ol{\sigma}^{\p\p}}^*\bigl(\O_{\P(({\cal B}^\vee_{j_0|U^{\p\p\p}}/{\cal B}^\vee_{j_0-1|U^{\p\p\p}})^\vee)}(1)\bigr)=
\O_{U^{\p\p\p}}\bigl(-(D_{|U^{\p\p}}+D^\p)\bigr).
$$
Now, let $\chi_*$ be the character of $H_{\ol{G}}(\la)$ as in Proposition~\ref{RamRamI}. By our strengthening of
the rationality properties, $\chi_*$ comes from a character of $H_{G}(\la)$ which we denote again by
$\chi_*$. We may view $\chi_*$ also as a character of $Q_{\GL(V)}(\la)$. 
The given $Q_{\GL(V)}(\la)$-bundle ${\cal Q}_{|U^{\p\p\p}}\subset {\cal P}_{|U^{\p\p\p}}$ and $\chi_*^{-1}$ define
a line bundle $\L_{\chi_*^{-1}}$, and ${\cal Q}_{|U^{\p\p\p}}$ and $W_{j_0}/W_{j_0-1}\otimes\C_{\chi_*^{-1}}$ define a vector bundle
$$
\widetilde{\cal B}^\vee\ \cong\ 
\bigl({\cal B}^\vee_{j_0|U^{\p\p\p}}/{\cal B}^\vee_{j_0-1|U^{\p\p\p}}\bigr)\otimes \L_{\chi_*^{-1}}
$$ over $U^{\p\p\p}$, so that 
$$
{\ol{\sigma}^{\p\p}}^*\bigl(\O_{\P(\widetilde{\cal B})}(1)\bigr)\ \cong\
\O_{U^{\p\p\p}}\bigl(-(D_{|U^{\p\p}}+D^\p)\bigr)\otimes \L^\vee_{\chi^{-1}_*}.
$$
Now, Proposition~\ref{Mummy} grants that the assumptions of Proposition~\ref{RamRamII} are satisfied,
so that we conclude that 
$$
0\ \ge\ \deg\Bigl(\O_{U^{\p\p\p}}\bigl(-(D_{|U^{\p\p}}+D^\p)\bigr)\Bigr)\ \ge\ \deg(\L_{\chi_*^{-1}}).
$$
By construction, $\chi_*^{-1}$ is just a \sl positive \rm multiple of $l_T(\la)$, so that also
$\deg(\L_{l_T(\la)})\le 0$. 
Set
$$
(\gamma^\p_1,...,\gamma^\p_{s+1})\q:=\q
\sum_{i=1}^s\alpha^\p_i\cdot\bigl(\underbrace{\rk(\A^\p_i)-r,...,\rk(\A^\p_i)-r}_{\rk(\A^\p_i)\times},
\underbrace{\rk(\A^\p_i),...,\rk(\A^\p_i)}_{(r-\rk(\A^\p_i))\times}\bigr).
$$ 
Using Example~\ref{CharComp}, we find
\begin{eqnarray*}
\deg(\L_{l_T(\la)}) &=& \sum_{i=1}^{s+1}\gamma^\p_i\Bigl(\deg\bigl(\det(\A^\p_{i})\bigr)-\deg\bigl(\det(\A^\p_{i-1})\bigr)\Bigr)
\\
&\stackrel{\deg(\A_{s+1})=0}{=} &
                        -\sum_{i=1}^{s}(\gamma^\p_{i+1}-\gamma^\p_i)\cdot\deg\bigl(\det(\A^\p_{i})\bigr)
\\
                    &=& -\sum_{i=1}^{s}\alpha^\p_i\cdot r \cdot\deg(\A^\p_{i})
\\
                    &\stackrel{(\ref{degreeequality})}=& -\sum_{i=1}^{s}\alpha_{s+1-i}\cdot r \cdot\deg(\A_{s+1-i})
\\
                    &=& -\sum_{i=1}^{s}\alpha_{i}\cdot r \cdot\deg(\A_{i}).		
\end{eqnarray*}
The last expression is the coefficient of the monomial of degree ${\dim(X)-1}$ in $M(\A^\bullet,\ul{\alpha})$, 
so that the claim
and hence the theorem is settled.
\end{proof}
\section{Construction of the Representation $\kappa$ and Applications}
Let $(\A,\tau)$ be an honest singular $G$-bundle. It will be our main problem to characterize among all weighted
filtrations $(\A^\bullet,\ul{\alpha})$ of $\A$ those which are associated with a reduction of $(\A,\tau)$ to a one
parameter subgroup of $G$.
\subsection{Definition of $\kappa$ and Elementary Properties}
\label{DefOfKap}
First, there exist a representation $\widetilde{\kappa}\colon \SL(V)\lra \GL(\widetilde{W})$ 
and an $\SL(V)$-equivariant embedding
$$
\eta\colon\  \SL(V)/G\hookrightarrow \widetilde{W},
$$
by \cite{Bo}, 1.12 Proposition. 
The representation $\widetilde{\kappa}$ may be extended to a representation
$\kappa\colon \GL(V)\lra \GL(W)$. By twisting the direct summands of $\kappa$
by suitable powers of the determinant --- which does not alter $\widetilde{\kappa}$, we may assume that all weights
of the action of $\C^*=\C^*\cdot \id_V$ be positive. Denote the resulting representation
again by $\kappa$.
The important features of this construction are summarized
in the following lemma. For $g\in\GL(V)$, $[g]$ stands for the image of $g\cdot G$ in $W$.
\begin{Lem} 
\label{FiniteIndex}
{\rm i)} The group $G$ has finite index in $\widetilde{G}$, the $\GL(V)$-stabilizer of $[e]$.
\par
{\rm ii)} For every $g\in \GL(V)$, the point $[g]$ is $\SL(V)$-polystable.
\par
{\rm iii)} A point $k\in \ol{\GL(V)\cdot [e]}\setminus \GL(V)\cdot [e]$ is not $\SL(V)$-semistable.
\end{Lem}
\begin{proof}
Ad i). We look at the isogeny $\SL(V)\times\C^*\lra \GL(V)$, $(h,z)\lma h\cdot z$. Note that the
$\SL(V)$-orbit of $[e]$ is, by construction, closed in $W$, and, in particular, does not contain
the origin. Since we assume that all weights of the $\C^*$-action be positive, the $\SL(V)$-orbit of $[e]$
intersects the $\C^*$-orbit of $[e]$ only in finitely many points. This implies that $G$ has finite index in the 
$(\SL(V)\times\C^*)$-stabilizer of $[e]$ and settles our claim.
\par
Ad ii). Since the $\SL(V)$-orbit of $[e]$ is closed in $W$, the commutative diagram
$$
\begin{CD}
\GL(V) \times W @>>> W
\\
@V (g\cdot .\cdot g^{-1}, V g\cdot)V @VV g\cdot V
\\
\GL(V)\times W @>>> W
\end{CD}
$$
shows that the $(g\cdot \SL(V)\cdot g^{-1})$-orbit of $[g]$ is closed and does not contain the origin. Finally,
$\SL(V)$ is normal in $\GL(V)$.
\par
Ad iii). We write $\ol{\GL(V)\cdot [e]}=\Spec(A)$. Since all the points $[g]$, $g\in \GL(V)$, are $\SL(V)$-polystable,
we get a dominant morphism
$$
\C^*\cong\GL(V)/\SL(V)\lra \Spec(A^{\SL(V)}).
$$
Thus, $A^{\SL(V)}\subset \C[f]$ with $f\colon \GL(V)/\SL(V)\lra\C$, $[g]\lma \det(g)$
(the function $f^{-1}$ is not regular at $0\in \ol{\GL(V)\cdot [e]}$). 
This shows that all elements in the closed
set
$\ol{\GL(V)\cdot [e]}\setminus \GL(V)\cdot [e]$ are nullforms.
\end{proof}
\begin{Rem}
These are the properties --- alluded to in the introduction --- which our construction shares, for example, with
the case of an adjoint group and the embedding $\GL(\frak g)/{\rm Aut}(\frak g)\hookrightarrow
\Hom(\frak g\otimes\frak g,\frak g)$. Our theory relies on these properties.
\end{Rem}
If we are given a point $[g]\in W$ and a one parameter subgroup $\la\colon \C^*\lra g\cdot
\widetilde{G}\cdot g^{-1}$, then $[g]$ is a fixed point for the induced $\C^*$-action and we have
$\mu(\la,[g])=0$. We will first establish a kind of converse to this trivial observation and then extend it
to weighted filtrations, thereby arriving at the necessary characterization of weighted filtrations 
$(\A^\bullet,\ul{\alpha})$ arising
from reductions to one parameter subgroups as those for which $\mu(\A^\bullet,\ul{\alpha};\phi)=0$.
Here, $\phi$ is as in Example \ref{weightcomputation}.
\par
Suppose we are given a point $[g]\in W$ and a one parameter subgroup
$\la\colon\C^*\lra\SL(V)$ with $\mu(\la,[g])=0$. Since $[g]$ is $\SL(V)$-polystable, it follows that
$\lim_{z\rightarrow \infty}\la(z)\cdot [g]=[g^\p]$ for some $[g^\p]\in \GL(V)$, and $\la$ is a one parameter
subgroup of the component of the identity of the $\GL(V)$-stabilizer of $[g^\p]$, that is, of 
$G^{g^\p}:=g^\p\cdot G\cdot{g^\p}^{-1}$. Therefore, $Q_{\GL(V)}(\la)\cap G^{g^\p}$ is the parabolic
subgroup $Q_{G^{g^\p}}(\la)$. By Example \ref{HomSpace}, there
is an element $u\in {\cal R}_u(Q_{\GL(V)}(\la))$,
such that $u\cdot [g^\p]= [g]$. As we have seen before, $Q_{\GL(V)}(\la)\cap G^{g^\p}$ 
is the parabolic
subgroup $Q_{G^{g^\p}}(\la)$. Observing that 
$u\cdot Q_{\GL(V)}(\la)\cdot u^{-1}=Q_{\GL(V)}(\la)$,
we conclude 
\begin{Prop}
\label{intersect}
If we are given a one parameter subgroup $\la\colon \C^*\lra \SL(V)$ and a point $[g]\in W$ with
$\mu(\la,[g])=0$, then
$$
Q_{\GL(V)}(\la)\cap g\cdot G\cdot g^{-1}=
 Q_{g\cdot G\cdot g^{-1}}(\la^\p)
$$
for some one parameter subgroup $\la^\p$ ($=u\cdot\la\cdot u^{-1}$ in the above notation)
of $g\cdot G\cdot g^{-1}$ with 
$Q_{\GL(V)}(\la^\p)=Q_{\GL(V)}(\la)$. More precisely, $\la$ and $\la^\p$ define the same weighted flag
in $V$.
\end{Prop}
\subsection{Characterization of Certain Weighted Filtrations}
\label{CharWeightFilt}
To simplify our arguments, we assume that 
$W=V_{\ul{a},\ul{b},\ul{c}}=V_{a_1,b_1,c_1}\oplus\cdots\oplus V_{a_n,b_n,c_n}$ 
(compare Example \ref{weightcomputation}).
Let $(\A,\tau)$ be an honest singular $G$-bundle, and let $\kappa\colon \GL(V)\lra \GL(W)$ be the representation 
constructed in the last section. Over the open set $U$ where $\A$ is locally free, we have
the reduction
$$
\sigma\colon\ U\lra \ul{\rm Isom}(V\otimes\O_U,\A^\vee_{|U})/G
$$
(cf.\ \cite{Schmitt}, Remark 3.3). Recall that ${\cal P}:=\ul{\rm Isom}(V\otimes\O_U,\A^\vee_{|U})$
is the principal $\GL(V)$-bundle corresponding to $\A^\vee_{|U}$ and that $\sigma$ and ${\cal P}\lra
{\cal P}/G$ define over $U$ the principal bundle ${\cal P}(\A,\tau)$, i.e., $\sigma$ encodes the rational
principal $G$-bundle associated with the singular principal $G$-bundle $(\A,\tau)$.
Using the representation $\kappa$, we get yet another description, namely $\sigma$ gives rise to a section
$$
\phi_{U}\colon\ U\lra \bigl((\A_{|U}^\vee)^{\otimes a_1}\bigr)^{\oplus b_1}\oplus\cdots\oplus 
\bigl((\A_{|U}^\vee)^{\otimes a_n}\bigr)^{\oplus b_n}\cong {\cal P}(\A,\tau)_{\kappa}
$$
which extends to a decoration
$$
\phi\colon\ \bigl(\A^{\otimes a_1}\bigr)^{\oplus b_1}\oplus\cdots\oplus 
\bigl(\A^{\otimes a_n}\bigr)^{\oplus b_n}\lra \O_X
$$
of type $(\ul{a},\ul{b},\ul{c})$.
Let $\beta$ be a reduction of $(\A,\tau)$ to the one parameter
subgroup $\la$ of $G$ and $(\A^\bullet_\beta,\ul{\alpha}_\beta)$ the associated weighted filtration of $\A$.
Our first contention is
\begin{Lem}
\label{CharWeightFilt1}
$$
\mu\bigl(\A^\bullet_\beta,\ul{\alpha}_\beta; \phi)\q=\q 0.
$$
\end{Lem}
\begin{proof}
If we apply the construction described in Example~\ref{weightcomputation}, the resulting section
$$
\phi^\p_{|\widetilde{U}}\colon\ \widetilde{U}\lra {\cal P}_{\kappa|\widetilde{U}}\cong
W\times\widetilde{U}\lra W,
$$
is just $x\lma [e]$, $x\in U$, so that the assertion is an obvious consequence of Formula~(\ref{weightcomp2}).
\end{proof}
Next, we come to the converse, i.e., to
\begin{Prop}
\label{CharWeightFilt2}
If $(\A^\bullet,\ul{\alpha})$ is a weighted filtration of $\A$ with 
$\mu\bigl(\A^\bullet,\ul{\alpha}; \phi)=0$, there exists a reduction $\beta$ to a one
parameter subgroup $\la\colon\C^*\lra G$ with
$$
(\A^\bullet,\ul{\alpha})\q=\q(\A^\bullet_\beta,\ul{\alpha}_\beta).
$$
\end{Prop}
\begin{proof}
Let $\eta$ be the generic point of $X$, $K$ the residue field at $\eta$, and ${\Bbb K}$ an algebraic closure
of $K$. Over ${\Bbb K}$, we may realize a situation as in Proposition \ref{intersect}. 
For this, let ${\cal P}:={\cal P}(\A,\tau)$ be the associated principal bundle. 
By applying the inverse 
of the procedure described in the introduction, the weighted filtration $(\A^\bullet, \ul{\alpha})$ defines,
over the big open subset $U^\p\subset U$ where all the $\A_i$ are subbundles, $i=1,...,s$, a weighted
filtration $(\A^{\p^\bullet}_{U^\p}, \ul{\alpha}^\p)$ of $\A^\vee_{|U^\p}$ by subbundles, i.e.,
$U^\p$ is the maximal open subset where $\A$, all the $\A_i$ and all the quotients $\A/\A_i$, $i=1,...,s$, 
are locally free. 
Then, for a suitable
parabolic subgroup $Q\subset\GL(V)$, this filtration defines a section 
$\beta\colon U^\p\lra \ul{\rm Isom}(V\otimes\O_{U^\p},\A^\vee_{|U^\p})/Q$ and, thus, a $Q$-bundle ${\cal Q}\subset
\ul{\rm Isom}(V\otimes\O_{U^\p},\A^\vee_{|U^\p})$.
\par 
The group $G$ acts on
itself, on $\widetilde{G}$, on $\SL(V)$, and on $\GL(V)$ by conjugation which provides us with the group schemes 
${\cal G}({\cal P})$, $\widetilde{\cal G}({\cal P})$, ${\cal SL}({\cal P})$, and ${\cal GL}({\cal P})$ over $U$. 
Over $\Spec({\Bbb K})$, this yields the groups $\rho_{\Bbb K}\colon G_{\Bbb K}\subset \SL_{\Bbb K}$,
with $G_{\Bbb K}:= {\cal G}({\cal P})\times_U \Spec({\Bbb K})$ and 
$\SL_{\Bbb K}:={\cal SL}({\cal P})\times_U \Spec({\Bbb K})$,
as well as $\widetilde{\rho}_{\Bbb K}\colon \widetilde{G}_{\Bbb K}\subset \GL_{\Bbb K}$,
with $\widetilde{G}_{\Bbb K}:= \widetilde{\cal G}({\cal P})\times_U \Spec({\Bbb K})$ and 
$\GL_{\Bbb K}:={\cal GL}({\cal P})\times_U \Spec({\Bbb K})$. 
If we choose, over a suitable finite extension
$\widetilde{K}$ of $K$, a trivialization ${\cal P}\times_U\Spec(\widetilde{K})\cong G\times_{\Spec(\C)}
\Spec(\widetilde{K})$,
we get induced trivializations $\widetilde{G}_{\Bbb K}\cong \widetilde{G}\times_{\Spec(\C)} \Spec({\Bbb K})$
and ${\GL}_{\Bbb K}\cong \GL(V)\times_{\Spec(\C)} \Spec({\Bbb K})$, such that $\widetilde{\rho}_{\Bbb K}
= \widetilde{\rho}\times_{\Spec(\C)} {\id}_{\Spec({\Bbb K})}$, $\widetilde{\rho}\colon \widetilde{G}\subset \GL(V)$.
Likewise, we get an induced representation $\kappa_{\Bbb K}\colon \GL_{\Bbb K}\lra \GL(W_{\Bbb K})$,
$W_{\Bbb K}:=W\otimes_\C {\Bbb K}\cong {\cal P}_{\kappa}\times_U \Spec({\Bbb K})$.
Under the identification of ${\cal P}_{\kappa}\times_U \Spec({\Bbb K})$ with $W_{\Bbb K}$,
$\phi_U$ and $\Spec({\Bbb K})\lra \eta\in U$ define the point $x_{\Bbb K}\in W_{\Bbb K}$. 
Then, $\widetilde{G}_{\Bbb K}$ is the $\GL_{\Bbb K}$-stabilizer
of $x_{\Bbb K}$, and the point $x_{\Bbb K}$ is $\SL_{\Bbb K}$-polystable. Next, the weighted filtration
$(\A^{\p^\bullet}_{U^\p}, \ul{\alpha}^\p)$ produces a weighted filtration $(\ol{V}^\bullet, \ul{\alpha}^\p)$
in $V_{\Bbb K}:=V\otimes_\C{\Bbb K}$. This weighted filtration can be obtained from a one
parameter subgroup $\la_{\Bbb K}\colon {\Bbb G}_m({\Bbb K})\lra \SL_{\Bbb K}$, and we have
$$
\mu_{\kappa_{\Bbb K}}\bigl(\la_{\Bbb K}, x_{\Bbb K}\bigr)\q=\q 0,
$$
by assumption.
Proposition~\ref{intersect} may now be applied. This means that there is a one parameter subgroup
$\widetilde{\la}_{\Bbb K}\colon {\Bbb G}_m({\Bbb K})\lra \widetilde{G}_{\Bbb K}$, such that
$$
Q_{\Bbb K}=Q_{\GL_{\Bbb K}}(\widetilde{\la}_{\Bbb K}).
$$
Here, $Q_{\Bbb K}:={\cal Q}\times_{U^\p} \Spec(\Bbb K)$. Using the $\C$-rational structure of $G_{\Bbb K}$,
we find a maximal torus $\ol{T}_{G_{\Bbb K}}$ in $G_{\Bbb K}$ of the form $T_G\times_{\Spec(\C)} \Spec({\Bbb K})$,
$T_G\subset G$ a maximal torus. Thus, we may find an element $g\in G_{\Bbb K}$ and a one parameter
$\widetilde{\la}\colon \C^*\lra G$, such that
$$
g\cdot Q_{G_{\Bbb K}}(\widetilde{\la}_{\Bbb K})\cdot g^{-1}\q =\q
Q_G(\widetilde{\la})\times_{\Spec(\C)} \Spec({\Bbb K}),
$$
whence
$$
g\cdot Q_{\Bbb K}\cdot g^{-1}\q =\q Q_{\GL(V)}(\widetilde{\la})\times_{\Spec(\C)} \Spec({\Bbb K}).
$$
This shows already that we may assume $Q=Q_{\GL(V)}(\widetilde{\la})$.
Since everything is defined over a finite extension of $K$, we have arrived at the following
\begin{Conc}
There are a finite Galois extension $K^\p/K$ and a trivialization $t\colon {\cal P}\times_U\Spec(K^\p)
\cong G\times_{\Spec(\C)}\Spec(K^\p)$, such that the section
\begin{eqnarray*}
\Spec(K^\p)&\stackrel{\beta\times_{U^\p}{\id}_{\Spec(K^\p)}}{\lra}&
\bigl(\ul{\rm Isom}(V\otimes \O_{U^\p}, \A^\vee_{|U^\p})/Q_{\GL(V)}(\widetilde{\la})\bigr)\times_{\Spec(\C)}\Spec(K^\p)
\\
&\cong& \bigl(\GL(V)/Q_{\GL(V)}(\widetilde{\la})\bigr)\times_{\Spec(\C)} \Spec(K^\p)
\end{eqnarray*}
lies in $(G/Q_G(\widetilde{\la}))\times_{\Spec(\C)} \Spec(K^\p)$.
\end{Conc}
Using the embedding $G/Q_G(\widetilde{\la})\hookrightarrow \GL(V)/Q_{\GL(V)}(\widetilde{\la})$, we find an embedding
$$
{\cal P}/Q_G(\widetilde{\la})\hookrightarrow \ul{\rm Isom}(V\otimes \O_{U^\p}, 
\A^\vee_{|U^\p})/Q_{\GL(V)}(\widetilde{\la}).
$$
We have just seen that $\beta\times_{U^\p}{\id}_{\Spec(K^\p)}$ lies in
$({\cal P}\times_U\Spec(K^\p))/Q_G(\widetilde{\la})$. But since everything is defined over $K$, we also must have
that $\beta\times_{U^\p}{\id}_{\Spec(K)}$ lies in ${\cal P}/Q_G(\widetilde{\la})$. This, of course, shows that
$\beta$ factorizes over ${\cal P}_{|U^\p}/Q_G(\widetilde{\la})$, and that implies the claim of the proposition.
\end{proof}
\section{Analysis of Semistability}
\label{AnStab}
\subsection{Slope Semistability and Mumford Semistability of $\A$}
Let $\A$ be a torsion free sheaf on $X$ and $(\A^\bullet,\ul{\alpha})$ a weighted filtration of $\A$.
Then, we set
\begin{eqnarray*}
L(\A^\bullet,\ul{\alpha})&=& \sum_{i=1}^s \alpha_i\bigl(\rk(\A_i)\cdot\deg(\A)-\rk(\A)\cdot\deg(\A_i)\bigr) 
\\
&=& \hbox{coefficient of $x^{\dim(X)-1}$ in $M(\A^\bullet,\ul{\alpha})$}.
\end{eqnarray*}
An honest singular $G$-bundle $({\cal A},\tau)$ is said to be \it slope-(semi)stable\rm,
if for every one parameter subgroup $\la\colon \C^*\lra G$ and every reduction $\beta$ of $({\cal A},\tau)$
to $\la$, we have
$$
L({\cal A}^\bullet_\beta, \ul{\alpha}_\beta)\q (\ge)\q 0.
$$
Then, one has the following implications
$$
\begin{array}{rcl}
\cal ({\cal A},\tau) \hbox{ is slope-stable} &\Longrightarrow& ({\cal A},\tau) \hbox{ is stable}
\\
&\Longrightarrow& ({\cal A},\tau) \hbox{ is semistable}
\\
&\Longrightarrow& ({\cal A},\tau) \hbox{ is slope-semistable}.
\end{array}
$$
\begin{Prop}
\label{equality99}
If the honest singular $G$-bundle $(\A,\tau)$ is slope-semistable, then the sheaf $\A$ is Mumford-semistable.
\end{Prop}
\begin{proof}
This is proved like Proposition~3.13 
in \cite{RamRam}. But, since we use other conventions and \cite{RamRam} contains some
sign errors, we will review the proof here.
Let ${\cal B}\subsetneq\A$ be any non-trivial saturated subsheaf. Set ${\cal B^\p}:=\ker(\A^\vee\lra {\cal B}^\vee)$.
We have to show that $\deg({\cal B})\le 0$. Over $U$, we have the principal $\GL(V)$-bundle
${\cal R}:=\ul{\rm Isom}(V\otimes\O_U, \A_{|U}^\vee)$, and, over the big open subset $U^\p\subset U$ 
where ${\cal B}^\p$ is a subbundle, it defines a section
$$
\sigma\colon\ U^\p\lra {\cal R}(\frak G, \alpha).
$$
Here, ${\frak G}:={\rm Gr}(V, i)$ is the Gra\ss mannian
of $i$-dimensional quotients of $V$ and $\alpha\colon\GL(V)\times {\frak G}
\lra {\frak G}$ is the usual action. The determinant $\L$ of the universal quotient is very ample.
We have a canonical linearization $\ol{\alpha}\colon \GL(V)\times\L\lra\L$. 
Note that the pullback of the line bundle ${\cal R}(\L,\ol{\alpha})\lra {\cal R}(\frak G,\alpha)$
is $\det(\A^\vee_{|U^\p}/{\cal B}_{|U^\p}^\p)$ which has degree $-\deg({\cal B}_{|U^\p}^\p)=-\deg({\cal B})$. 
If $\sigma(\eta)$ is semistable, then, by Proposition~\ref{RamRamII}, $-\deg({\cal B})\ge 0$, and we are
done.
\par
Otherwise, set $W^\vee:=H^0(\frak G,\L)$. 
As in the proof of Theorem~\ref{SemStabRed}, we find a one parameter subgroup $\widetilde{\la}\colon
\C^*\lra G$, a subquotient $W^{\p}$ of $W$ which inherits an $L_G(\widetilde{\la})$-module structure, and a section
$$
\sigma^\p\colon U^{\p\p}\lra {\cal R}_{|U^{\p\p}}\bigl(\P(W^{\p^\vee}),\alpha^\p\bigr),
$$
such that the pullback of the line bundle ${\cal R}_{|U^{\p\p}}\bigl(\O_{\P(W^{\p^\vee})}(1),\ol{\alpha}^\p\bigr)
\lra {\cal R}_{|U^{\p\p}}\bigl(\P(W^{\p^\vee}),\alpha^\p\bigr)$ is $\det({\cal B}^\p_{|U^{\p\p}})^\vee(-D)$
for some effective divisor $D$. Here, $\alpha^\p\colon L_{G}(\widetilde{\la})\times \P(W^{\p^\vee})\lra \P(W^{\p^\vee})$
and $\ol{\alpha}^\p\colon L_{G}(\widetilde{\la})\times \O_{\P(W^{\p^\vee})}(1)\lra 
\O_{\P(W^{\p^\vee})}(1)$ are induced by the $L_G(\widetilde{\la})$-module structure of $W^\p$.
Moreover, there is a character $\chi_*$ of $L_G(\widetilde{\la})$, such that $\sigma^\p(\eta)$ is semistable
w.r.t.\ the given linearization twisted by $\chi_*$. As before, one finds
$$
-\deg({\cal B})\ge \deg(\det({\cal B}^\p_{|U^{\p\p}})^\vee(-D))\ge \deg(\L_{\chi^{-1}_*}).
$$
Again, one checks that $\deg(\L_{\chi^{-1}_*})$ is a positive multiple of $L(A^\bullet_\beta,\ul{\alpha}_\beta)$
for some reduction $\beta$ to a one parameter subgroup of $G$. Our assumption, therefore, implies
that $\deg(\L_{\chi^{-1}_*})$ is non-negative, and we are done.
\end{proof}
\subsection{Comparison with Semistability of the Associated Decorated Sheaf}
First, we recall 
from our paper \cite{Schmitt} the definition of $\eps$-(semi)stability of a pair $(\A,\tau)$ where
$\eps\in\Q[x]$ is a positive polynomial of degree at most $\dim(X)-1$. For any $s>0$, first define
for every  $\ul{d}=(d_1, \cdots, d_s)$ with $\sum id_i = s!$ 
$$
\phi_{\ul{d}}\colon \bigotimes_{i=1}^{s} (V \otimes A)^{\otimes d_i i} \lra
\bigotimes_{i=1}^{s} \ {\rm Sym}^{d_i} 
\bigl({\rm Sym}^i  (V \otimes A)\bigr) \lra
\bigotimes_{i=1}^{s} {\rm Sym}^{d_i} ({\rm Sym}^i 
(V \otimes A))^G \lra {\cal O}_X
$$
where the first one is the standard map and the second one the Reynolds operator.
These add to
$$
\phi^s_\tau\colon \bigl((V \otimes A)^{\otimes s!}\bigr)^{\oplus N} 
\stackrel{\sum_{\ul{d}}\phi_{\ul{d}}}{\lra}{\cal O}_X.
$$
For $s\gg 0$, $\phi^s_\tau$
will be a non-zero homomorphism which defines $\tau$ uniquely. For any weighted filtration $(\A^\bullet,
\ul{\alpha})$, we define
$$
\mu(\A^\bullet,\ul{\alpha},\tau):=\frac{1}{s!}\cdot \mu(\A^\bullet,\ul{\alpha}; \phi^s_\tau).
$$
It is easy to see that this is well defined. Then, $(\A,\tau)$ is called \it $\eps$-(semi)stable\rm,
if for every weighted filtration $(\A^\bullet, \ul{\alpha})$ of $\A$, the inequality
$$
M(\A^\bullet,\ul{\alpha})+\eps\cdot\mu(\A^\bullet,\ul{\alpha},\tau)(\succeq)0
$$
holds.
\begin{Rem}
\label{SchB}
If $(\A,\tau)$ is $\eps$-semistable, then one has
$$
\mu({\cal B})\le \mu(\A)+\eps_0\cdot (r-1),
$$
for every saturated subsheaf $0\subsetneq {\cal B} \subsetneq \A$ (\cite{Schmitt}, (3.20)).
\end{Rem}
As we have described in Section~\ref{CharWeightFilt}, an honest singular principal $G$-bundle $(\A,\tau)$
defines an associated torsion free sheaf $(\A,\phi)$ with a decoration of type $(\ul{a},\ul{b},\ul{c})$.
Some important properties of the semistability concept are summarized in 
\begin{Prop}
\label{semstab100}
Let $\eps\in\Q[x]$ be a positive polynomial of degree exactly $\dim(X)-1$. Then the following properties
hold true.
\begin{enumerate}
\item[\rm 1.] An honest singular $G$-bundle $(\A,\tau)$ is (semi)stable, if and only if the associated
              decorated sheaf $(\A,\phi)$ is $\eps$-(semi)stable (as defined in Section {\rm\ref{DecSh}}).
\item[\rm 2.] If $(\A,\tau)$ is a (semi)stable honest singular $G$-bundle, then it is
              $\eps$-(semi)stable (in the sense of {\rm\cite{Schmitt}}).
\end{enumerate}
\end{Prop}
\begin{proof} 
Ad 1. By the Lemma~\ref{CharWeightFilt1}, the 
$\eps$-(semi)stability of $(\A,\phi)$ clearly implies the (semi)stability of $(\A,\tau)$. 
If $(\A,\tau)$ is (semi)stable, Proposition~\ref{equality99} shows that $\A$ is a Mumford-semistable torsion free sheaf.
Therefore, for every weighted filtration $(\A^\bullet,\ul{\alpha})$, the number $L(\A^\bullet,\ul{\alpha})$
is non-negative. Thus, if we have $\mu(\A^\bullet,\ul{\alpha};\phi)>0$, then
$M(A^\bullet,\ul{\alpha})+\eps\cdot \nu(\A^\bullet,\ul{\alpha};\phi)\succ 0$. On the other hand, if 
$\mu(\A^\bullet,\ul{\alpha};\phi)=0$, then we may apply Proposition~\ref{CharWeightFilt2} in order to see that
the condition $M(A^\bullet,\ul{\alpha})(\succeq) 0$ follows from the (semi)stability of $(\A,\tau)$.
\par
Ad 2. We use the rational map 
$$
\P\bigl(\Hom(V,V)\bigr)
\dasharrow
\P\bigl(\Hom(V,V)\bigr)\catqot G
\lra 
\P\bigl(\Hom(V,V)\bigr)\catqot\widetilde{G}
\dasharrow 
V_{\ul{a},\ul{b},\ul{c}}\catqot\C^*
$$
which is defined on $\P({\rm Isom}(V,V))$, the last map coming from
$\GL(V)/\widetilde{G}\subset V_{\ul{a},\ul{b},\ul{c}}\setminus\{0\}$. It extends over a big open subset
$U\subset \P\bigl(\Hom(V,V)\bigr)$, and the pullback of any ample line bundle $\L$ on the weighted projective
space $V_{\ul{a},\ul{b},\ul{c}}\catqot\C^*$ extends to an $(\SL(V)\times \widetilde{G})$-linearized
ample line bundle $\O_{\P(\Hom(V,V))}(l)$. The $\SL(V)$-invariant global sections of $\L$ define 
$(\SL(V)\times \widetilde{G})$-invariant global sections of $\O_{\P(\Hom(V,V))}(l)$ and, thus,
$G$-invariant global sections of $\n^{\otimes l^\p}$, with $\n$ the polarization on
$\P\bigl(\Hom(V,V)\bigr)\catqot G$ induced by the $G$-linearized ample line bundle 
$\O_{\P(\Hom(V,V))}(1)$. This shows that there is a fixed positive rational number $\eta$, such that
$$
\mu_\n(\la, [t])\q\ge \q \eta\cdot \mu_{\L}(\la, f([t])),\q \forall\ \la\colon\C^*\lra\SL(V),\ 
t\in {\rm Isom}(V,V),
$$
$f\colon \P\bigl(\Hom(V,V)\bigr)\catqot G
\dasharrow V_{\ul{a},\ul{b},\ul{c}}\catqot\C^*$ being the rational map defined before.
By the same token, we may find a positive rational number $\eta^\p$, such that
$$
M(\A^\bullet,\ul{\alpha})+\eps\cdot \mu(\A^\bullet,\ul{\alpha};\tau)
\q\succeq\q M(\A^\bullet,\ul{\alpha})+\eta^\p\cdot\eps\cdot \nu(\A^\bullet,\ul{\alpha};\phi).
$$
Since $(\A,\phi)$ is also $(\eta^\p\cdot\eps)$-(semi)stable, by 1., we are done.
\end{proof}
\subsection{Ramanathan's Definition of Semistability}
Here, we assume that $G$ be semisimple. A \it rational principal $G$-bundle over $X$ \rm is a pair
$(U,{\cal P})$ which consists of a big open subset $U\subset X$ and a principal $G$-bundle ${\cal P}$
over $U$. A rational principal $G$-bundle $(U,{\cal P})$ is called \it Ramanathan-(semi)stable\rm, if for
every big open subset $U^\p\subset U$, every parabolic subgroup $Q$,  every antidominant character $\chi$
of $Q$, and every reduction $\beta\colon U^\p\lra {\cal P}_{|U^\p}/Q$, one has $\deg(\L_\chi)\ge 0$.
Here, $\L_{\chi}$ is the line bundle on $U^\p$ associated with $\beta$ and the $G$-linearized (ample) line
bundle on $G/Q$ assigned to the $Q$-bundle $G\lra G/Q$ via $\chi$. Since every antidominant character
is induced --- by the procedure outlined in Example \ref{CharComp} --- by a one parameter
subgroup $\la\colon \C^*\lra G$, we see
\begin{Prop}
An honest singular principal $G$-bundle $(\A,\tau)$ is slope (semi)sta\-ble, if and only the associated rational
principal bundle $(U,{\cal P}(\A,\tau))$ is Ramanathan-(semi)stable.
\end{Prop} 
\begin{Rem}
For general reductive groups, the antidominant characters 
correspond to one parameter subgroups of $[G,G]$.
\end{Rem}
\section{Semistable Reduction}
In this section, we will present the proof of the main theorem as stated in the introduction. The strategy of
proof is analogous to the one used in \cite{Schmitt} for adjoint groups and the adjoint representation.
\subsection{Singular $\widetilde{G}$-bundles}
Let $\widetilde{G}$ be the group defined in Section~\ref{DefOfKap}
and $\widetilde{\rho}\colon \widetilde{G}\lra \GL(V)$ its natural inclusion into $\GL(V)$. Since $\widetilde{\rho}(\widetilde{G}^0)
=\rho(G)\subset \SL(V)$, Proposition 1.15 in Chapter 1.4 of \cite{GIT} grants that we may develop a theory 
of singular principal
$\widetilde{G}$-bundles as well. Again, we have the following alternative:
\begin{Lem}
\label{honest}
Let $(\A,\widetilde{\tau})$ be a singular principal $\widetilde{G}$-bundle, $U$ the maximal open
subset where $\A$ is locally free, and $\widetilde{\sigma}\colon U\lra \ul{\Hom}(V\otimes\O_U,\A_{|U}^\vee)\catqot
\widetilde{G}$ the section defined by $\widetilde{\tau}_{|U}$. Then, $\widetilde{\sigma}$ is entirely contained
in either $\ul{\rm Isom}(V\otimes\O_{U},\A_{|U}^\vee)/\widetilde{G}$ or its complement.
In the former case, we call 
$(\A,\widetilde{\tau})$ an {\rm honest singular principal $\widetilde{G}$-bundle}.
\end{Lem}
\begin{proof} Let $h\colon \ul{\Hom}(V\otimes\O_{U},\A_{|U}^\vee)\lra U$ be the $\Hom$-bundle. Then, there
is the universal homomorphism $f\colon V\otimes\O_{\ul{\Hom}(V\otimes\O_{U},\A_{|U}^\vee)}\lra h^*(\A_{|U}^\vee)$.
Choosing a trivialization for $\det(\A)$ and a basis for $V$, $\det(f)$ becomes an $\SL(V)$-invariant
global section of 
$\O_{\ul{\Hom}(V\otimes\O_{U},\A_{|U}^\vee)}$. Since the image of $\widetilde{G}$ under the determinant
homomorphism is a finite, whence cyclic, subgroup of $\C^*$, some power of this function is, in fact,
$\widetilde{G}$-invariant, and, thus, defines a function $\frak d\in \Gamma(\ul{\Hom}(V\otimes\O_{U},\A_{|U}^\vee)
\catqot G,\O_{\ul{\Hom}(V\otimes\O_{U},\A_{|U}^\vee)\catqot G})$. Then, $\widetilde{\sigma}^*(\frak d)$
is a constant function, because $U$ is a big open subset of $X$. This implies the claim. 
\end{proof}
Next, let $S$ be any scheme and $(\A_S,\tau_S)$ a family of singular principal $G$-bundles.
We obtain the associated family $(\A_S,\widetilde{\tau}_S)$ of singular principal $\widetilde{G}$-bundles
with
$$
\widetilde{\tau}_S\colon\ {\Sym}^*(V\otimes\A_S)^{\widetilde{G}}
\subset
{\Sym}^*(V\otimes\A_S)^{G}\stackrel{\tau_S}{\lra} \O_{S\times X}.
$$
\begin{Rem} Note that ${\Sym}^*(V\otimes\A_S)^{\widetilde{G}}$ is the algebra of invariants in
${\Sym}^*(V\otimes\A_S)^{G}$ for the induced action of the finite group $\widetilde{G}/G$.
\end{Rem}
\begin{Lem}
\label{equality513}
{\rm i)} A singular
principal $G$-bundle $(\A,\tau)$ is honest, if and only if the associated 
singular
$\widetilde{G}$-bundle $(\A,\widetilde{\tau})$ is honest.
\par
{\rm ii)} For every positive polynomial $\eps\in \Q[x]$ of degree at most $\dim(X)-1$, a singular
principal $G$-bundle $(\A,\tau)$ is $\eps$-(semi)stable, if and only if the associated 
singular
$\widetilde{G}$-bundle $(\A,\widetilde{\tau})$ is $\eps$-(semi)stable.
\end{Lem}
\begin{proof} 
The first assertion is obvious. Since $\P(\Hom(V, V))\catqot \widetilde{G}$ is the quotient
of $\P(\Hom(V,\allowbreak V))\catqot G$ by the finite group $\widetilde{G}/G$, Lemma~\ref{equality500} implies
$$
\mu(A^\bullet,\ul{\alpha};\widetilde{\tau})\ =\ \mu(\A^\bullet,\ul{\alpha}; \tau)
$$
for every weighted filtration $(\A^\bullet,\ul{\alpha})$ of $\A$.
\end{proof}
\subsection{The Proof of the Main Theorem}
For the rest, we use the representation $\kappa\colon \GL(V)\lra
\GL(W)$, with $W=V_{\ul{a},\ul{b},\ul{c}}$ for appropriate tuples $\ul{a}$, $\ul{b}$, and $\ul{c}$, 
as constructed in Section~\ref{DefOfKap}. Moreover, we fix a Hilbert polynomial $P$ and
a positive polynomial $\eps\in\Q[x]$ of degree exactly $\dim(X)-1$. Then, the set of torsion free
sheaves $\A$ with Hilbert polynomial $P(\A)=P$
for which there exists either an $\eps$-semistable singular principal $\widetilde{G}$-bundle
$(\A, \widetilde{\tau})$  or an $\eps$-semistable torsion free sheaf $(\A,\phi)$ 
with a decoration of type
$(\ul{a},\ul{b},\ul{c})$ is bounded. This follows from Maruyama's boundedness result (\cite{HL}, Theorem 3.3.7)
and Remarks \ref{MarB}, ii), and \ref{SchB}. 
Therefore, we may choose an integer $n\gg 0$, such that
$\A(n)$ is globally generated and $H^i(\A(n))=0$, $i>0$,  for every such sheaf $\A$,
and such that all the
necessary GIT-constructions may be performed relative to the quasi-projective quot-scheme $\frak Q^0$.
Here, $\frak Q^0$ parameterizes quotients
$q\colon Y\otimes\O_X(-n)\lra\A$, such that $\A$ is a torsion free sheaf with Hilbert polynomial $P$
and $H^0(q(n))\colon Y\lra H^0(\A(n))$ is an isomorphism, $Y$ being a complex vector space
of dimension $P(n)$. In the following, the quotient maps involved, such as $q_S\colon Y\otimes\pi_X^*\O_X(-n)
\lra\A_S$, will remain unaffected by the constructions which will be carried out, 
whence they won't be mentioned.
\par
Now, as in Section 5.1 of \cite{Schmitt}, we may construct the parameter spaces $\frak Y(G)$ and
$\frak Y(\widetilde{G})$ over $\frak Q^0$ for singular principal $G$-bundles and singular
principal $\widetilde{G}$-bundles, respectively. Let us briefly recall the construction for $G$.
For large $s$, we start with 
$$\overline{\frak Y}={\frak Q}^0 \times \bigoplus_{i=1}^{s}{\rm 
Hom} \Bigl({\rm Sym}^i (V \otimes Y), H^0 \bigl({\cal O}_X (in)\bigr)\Bigr).
$$
Note that, over $\overline{\frak Y}\times X$, 
there are universal homomorphisms
$$
\widetilde{\varphi}^i:{\rm Sym}^i (V \otimes Y) \otimes {\cal O}_
{\overline{\frak Y} \times X} \to H^0 ({\cal O}_X (in)) \otimes 
{\cal O}_{\overline{\frak Y} \times X},\q i=1,..., s.
$$
Let $\varphi^i = {\rm ev} \circ \tilde{\varphi}^i$ be the composition of
$\widetilde{\phi}^i$ with the evaluation map  
${\rm ev}\colon H^0(\O_X(in))\otimes\O_{\ol{\frak Y}\times X}\lra \pi_X^*{\cal O}_X (in)$, $i=1,...,s$.  
We twist $\varphi^i$ by $\id_{\pi_X^*\O_X(-in)}$ and put the resulting maps together to  the  homomorphism 
$$
\varphi\colon {\cal V}_{\overline{\frak Y}} := 
\bigoplus_{i=1}^s {\rm Sym}^i\bigl(V\otimes Y\otimes \pi_X^* {\cal
O}_X(- n)\bigr) \lra {\cal O}_{\overline{\frak Y} 
\times X}.
$$ 
Next, $\phi$ yields a homomorphism of ${\cal O}_{\overline{\frak Y}\times X}$-algebras
$$
\widetilde{\tau}_{\overline{\frak Y}}\colon {\rm Sym}^* \bigl({\cal V}_{\overline{\frak Y}}\bigr) \lra
{\cal O}_{\overline{\frak Y} \times X}.
$$
On the other hand, there is a surjective homomorphism
$$
\beta\colon {\rm Sym}^* ({\cal V}_{\overline{\frak Y}}) \lra {\rm Sym}^* (V \otimes 
\pi^*{\cal A}_{{\frak  Q}^0})^G
$$
of graded algebras, where the left hand algebra is graded by assigning the weight $i$
to the elements in ${\rm Sym}^i(V\otimes Y\otimes \pi_X^*\O_X(-n))$.  
Here, $\pi\colon \ol{\frak Y}\times X\lra {\frak Q}^0\times X$ is the natural projection, and $\A_{{\frak Q}^0}$
is the universal quotient on ${\frak Q}^0\times X$.
The parameter space ${\frak Y}(G)$ is
defined by the condition that $\widetilde{\tau}_{\overline{\frak Y}}$ factorize over
$\beta$, i.e., setting 
$
{\cal A}_{{\frak Y}(G)}:=(\pi^*{\cal A}_{{\frak  Q}^0})_{|{\frak Y}(G)\times X}
$,
there be a homomorphism
$$
{\tau}_{{\frak Y}(G)} \colon {\rm Sym}^* (V \otimes {\cal A}_{{\frak Y}(G)})^G \lra  
{\cal O}_{{\frak Y}(G)\times X}
$$ 
with $\widetilde{\tau}_{\overline{\frak Y}|{\frak Y}(G)\times X}={\tau}_{{\frak Y}(G)}\circ\beta$.
Formally, ${\frak Y}(G)$ is defined as the scheme theoretic intersection of the closed subschemes
$$
{\frak Y}^d:=\Bigl\{\,y\in \ol{\frak Y}\,|\, \widetilde{\tau}^d_{\overline{\frak Y}|\{y\}\times X}
\colon \ker\beta_{|\{y\}\times X}^d\lra \O_X\ \hbox{\rm is trivial}\,\Bigr\}, \q d\ge 0.
$$
The family
$(\A_{{\frak Y}(G)}, \tau_{{\frak Y}(G)})$ is the \it universal family of singular
principal $G$-bundles parameterized by ${\frak Y}(G)$\rm. Similarly, we obtain ${\frak Y}(\widetilde{G})$.
\par
The homomorphism
$$
{\rm Sym}^* (V \otimes {\cal A}_{{\frak Y}(G)})^{\widetilde {G}}
\subset {\rm Sym}^* (V \otimes {\cal A}_{{\frak Y}(G)})^G \stackrel{{\tau}_{{\frak Y}(G)}}{\lra} 
{\cal O}_{{\frak Y}(G)\times X}
$$ 
leads, by the universal
property of $\frak Y(\widetilde{G})$, to a $\GL(Y)$-equivariant morphism
$$
\begin{CD}
\frak Y(G) @> \frak s>> \frak Y(\widetilde{G})
\\
@VVV @VVV
\\
\frak Q^0 @= \frak Q^0.
\end{CD}
$$
There is the induced commutative diagram
$$
\begin{CD}
\frak Y(G)\catqot\C^* @> \ol{\frak s}>> \frak Y(\widetilde{G})\catqot\C^*
\\
@VVV @VVV
\\
\frak Q^0 @= \frak Q^0.
\end{CD}
$$
Now, $\frak Y(G)\catqot\C^*\lra \frak Q^0$ is a proper morphism, so that the last diagram implies that
$\ol{\frak s}$ is a proper morphism, too. There are the $\GL(Y)$-invariant open subschemes $\frak Y^{\rm h}(G)$
and $\frak Y^{\rm h}(\widetilde{G})$, parameterizing the honest singular $G$- and $\widetilde{G}$-bundles,
respectively. Note that $\frak s^{\rm h}:={\frak s}_{|\frak Y^{\rm h}(G)}$ is quasi-finite
and  that $\frak s^{-1}(\frak Y^{\rm h}(\widetilde{G}))=\frak Y^{\rm h}(G)$, by Lemma~\ref{equality513}, i).
Thus, we also get a finite morphism 
$$
\ol{\frak s}^{\rm h}\colon \frak Y^{\rm h}(G)\catqot\C^*
\lra \frak Y^{\rm h}(\widetilde{G})\catqot\C^*.
$$
Finally, let $(\A,\widetilde{\tau}_S)$ be a family of honest singular principal $\widetilde{G}$-bundles
parameterized by $S$. Let ${\cal U}_{\rm lf}\subset S\times X$ be the maximal open subset where $\A_S$ is locally
free. Note that ${\cal U}_{\rm lf}\cap \{s\}\times X$ is the maximal open subset where
$\A_{S|\{s\}\times X}$ is locally free (\cite{HL}, Lemma 2.1.7). Then, over ${\cal U}_{\rm lf}$, we have the section
$$
\widetilde{\sigma}\colon {\cal U}_{\rm lf}\lra \ul{\rm Isom}(V\otimes \O_{{\cal U}_{\rm lf}},
\A^\vee_{S|{\cal U}_{\rm lf}})/\widetilde{G}
$$
induced by $\widetilde{\tau}_{S|{\cal U}_{\rm lf}}$. As in Section \ref{CharWeightFilt}, 
using the representation $\kappa$, the section $\widetilde{\sigma}$
gives rise to a section
$$
\sigma^\p\colon {\cal U}_{\rm lf}\lra \A^\vee_{S|{\cal U}_{\rm lf}; \ul{a},\ul{b},\ul{c}}
$$
or, dually,
$$
\phi_{{\cal U}_{\rm lf}}\colon \A_{S|{\cal U}_{\rm lf};\ul{a},\ul{b},\ul{c}} \lra \O_{{\cal U}_{\rm lf}}.
$$
This yields
$$
\phi_S\colon \A_{S;\ul{a},\ul{b},\ul{c}}\lra j_*\Bigl(\A_{S|{\cal U}_{\rm lf};\ul{a},\ul{b},\ul{c}}\Bigr)
\stackrel{j_*(\phi_{{\cal U}_{\rm lf}})}{\lra}
j_*(\O_{{\cal U}_{\rm lf}})=\O_{S\times X},
$$
$j\colon {\cal U}_{\rm lf}\lra S\times X$ being the inclusion. According to Maruyama (\cite{Maru}, bottom of page 112)
the equality $j_*(\O_{{\cal U}_{\rm lf}})=\O_{S\times X}$ is seen as follows: By \cite{EGAIV2}, Theorem (5.10.5)
--- applied to $S\times X$ and $Z:=(S\times X)\setminus {\cal U}_{\rm lf}$ (which is stable under specialization
\cite{HL}, Lemma 2.1.7) ---, one has to show that $\inf_{x\in Z}{\rm depth}(\O_{S\times X, x})\ge 2$. Since $X$ is smooth, the morphism
$\pi_S\colon S\times X\lra X$ is smooth. Thus, by \cite{EGAIV4}, Proposition (17.5.8),
$$
\dim(\O_{S\times X, x})-{\rm depth}(\O_{S\times X, x})= \dim(\O_{S,s})-{\rm depth}(\O_{S,s}),
$$
for every point $x\in S\times X$ and $s:=\pi_S(x)$.
This implies
\begin{equation}
\label{depth}
{\rm depth}(\O_{S\times X, x})\ge \dim(\O_{S\times X, x})-\dim(\O_{S,s})=\dim(\O_{\pi_S^{-1}(s),x}).
\end{equation}
Since for any point 
$x\in \pi_S^{-1}(s)$, one has $\dim\O_{\pi_S^{-1}(s),x}={\rm codim}_{\pi_S^{-1}(s)}(\ol{\{x\}})$,
we derive ${\rm depth}(\O_{S\times X, x})\ge 2$ for every point $x\in Z$ from the fact that 
${\rm codim}_{\pi_S^{-1}(s)}(Z\cap \pi_S^{-1}(s))\ge 2$ and (\ref{depth}).
\par
Denote by ${\frak W}_{\ul{a},\ul{b},\ul{c}}\lra \frak Q^0$ the parameter space for torsion free sheaves 
$(\A,\phi)$ with a decoration of
type $(\ul{a},\ul{b},\ul{c})$. The above construction then leads to a $\GL(Y)$-equivariant injective morphism
$$
\begin{CD}
\frak Y^{\rm h}(\widetilde{G}) @> \frak t>> \frak W_{\ul{a},\ul{b},\ul{c}}
\\
@VVV @VVV
\\
\frak Q^0 @= \frak Q^0.
\end{CD}
$$
Note that --- on closed points --- the morphism $\frak t\circ\frak s^{\rm h}$ is just the assignment which associates 
to an honest
singular principal $G$-bundle $(\A,\tau)$ the corresponding torsion free sheaf $(\A,\phi)$ with a decoration of
type $(\ul{a},\ul{b},\ul{c})$ that had been used in Section \ref{CharWeightFilt}.
In the following, ${\frak W}^{\eps\rm-ss}_{\ul{a},\ul{b},\ul{c}}$ stands for the open subset which parameterizes 
the $\eps$-semistable
torsion free sheaves with a decoration of type $(\ul{a},\ul{b},\ul{c})$.
We will prove
\begin{Thm}
\label{MainAux}
The restricted morphism
$$
\frak t^\p:=\frak t_{|t^{-1}({\frak W}^{\eps\rm-ss}_{\ul{a},\ul{b},\ul{c}})}\colon 
\frak t^{-1}\Bigl({\frak W}^{\eps\rm-ss}_{\ul{a},\ul{b},\ul{c}}\Bigr)
\lra {\frak W}^{\eps\rm-ss}_{\ul{a},\ul{b},\ul{c}}
$$
is proper, whence finite.
\end{Thm}
\begin{Rem}
Note that, by Proposition \ref{semstab100} and Lemma \ref{equality513}, ii), 
$({\frak t\circ\frak s^{\rm h}})^{-1}({\frak W}^{\eps\rm-ss}_{\ul{a},\ul{b},\ul{c}})$
is contained in the locus of $\eps$-semistable singular $G$-bundles.
\end{Rem}
This theorem implies that $\frak t^{-1}({\frak W}^{\eps\rm-ss}_{\ul{a},\ul{b},\ul{c}})$ admits
a projective $\GL(Y)$-quotient, or, equivalently, that
$\bigl(\frak t^{-1}({\frak W}^{\eps\rm-ss}_{\ul{a},\ul{b},\ul{c}})\bigr)\catqot\C^*$
admits a projective $\SL(Y)$-quotient. Here, one uses Lemma 5.1 of \cite{Ramanathan}.
By the finiteness of $\ol{\frak s}^{\rm h}$ and the same lemma, this
implies that
$$
\ol{\frak s}^{\rm h^{-1}}
\biggl(
\Bigl(\frak t^{-1}\bigl({\frak W}^{\eps\rm-ss}_{\ul{a},\ul{b},\ul{c}}\bigr)\Bigr)\catqot\C^*
\biggr)
$$
admits a projective $\SL(Y)$-quotient, too. But, by Proposition \ref{semstab100}, the latter set is just 
$\bigl(\frak Y^{\rm h, ss}(G)\bigr)\catqot\C^*$, $\frak Y^{\rm h, ss}(G)$ being the open subset of semistable
honest singular $G$-bundles. Therefore, ${\frak Y}^{\rm h,ss}(G)$ admits a projective $\GL(Y)$-quotient,
and that concludes the proof of the main theorem.\qed
\begin{proof}[Proof of Theorem~\ref{MainAux}]
Let $\frak I\subset {\frak W}^{\eps\rm-ss}_{\ul{a},\ul{b},\ul{c}}$ 
be the scheme theoretic image of $\frak t^{-1}({\frak W}^{\eps\rm-ss}_{\ul{a},\ul{b},\ul{c}})$
under $\frak t$. Our first contention is
\begin{Claim}
The point set underlying $\frak I$ is just the set theoretic image of 
$\frak t^{-1}({\frak W}^{\eps\rm-ss}_{\ul{a},\ul{b},\ul{c}})$ under $\frak t$.
\end{Claim}
The claim is seen as follows.
Let $\frak J$ be the set theoretic image. This is a constructible subset of $\frak W_{\ul{a},\ul{b},\ul{c}}$.
As such, it contains an open subset $\frak U$ of its closure, which is the closed subset underlying
$\frak I$. Now, suppose $w\in \frak I$ is a closed point. By what we have just observed, we may find
a smooth affine curve $C$ with a fixed point $*\in C$ and a map $f\colon C\lra \frak I$, such that
$f(C\setminus\{*\})\subset {\frak U}\subset \frak J$ and $f(*)=w$. Let $(\A_C,\phi_C)$ be 
the pullback of the universal
family on $\frak W_{\ul{a},\ul{b},\ul{c}}\times X$ to $C\times X$ and ${\cal U}_{\rm lf}\subset C\times X$ 
the maximal open subset over which $\A_C$ is locally free. We then have a section
$$
\sigma^\p\colon {\cal U}_{\rm lf}\lra \A^\vee_{C|{\cal U}_{\rm lf}; \ul{a},\ul{b},\ul{c}}.
$$
The vector bundle $\A^\vee_{C|{\cal U}_{\rm lf}}$ and the $\GL(V)$-varieties
$\GL(V)/\widetilde{G}=\GL(V)\cdot [e] \hookrightarrow W$ and $\ol{GL(V)\cdot [e]}
\hookrightarrow W$ define a locally closed subscheme 
${\cal S}\subset \A^\vee_{C|{\cal U}_{\rm lf}; \ul{a},\ul{b},\ul{c}}$
and a closed subscheme $\ol{\cal S}\hookrightarrow \A^\vee_{C|{\cal U}_{\rm lf}; \ul{a},\ul{b},\ul{c}}$.
By construction, the image of ${\cal U}_{\rm lf}\cap (C\setminus\{*\})\times X$ under $\sigma^\p$
lies in ${\cal S}$, whence 
$\sigma^\p({\cal U}_{\rm lf})$ lies in $\ol{\cal S}$. Now, Lemma~\ref{FiniteIndex}, iii),
and Theorem~\ref{SemStabRed} imply that there is also a non-empty open subset $U^\p\subset \{*\}\times X$,
such that $\sigma^\p(U^\p)\subset {\cal S}$. Thus, we have found an open subset 
${\cal U}^\p\subset {\cal U}_{\rm lf}$ which is big in $C\times X$, such that 
$\sigma:=\sigma^\p_{|{\cal U}^\p}$ factorizes over ${\cal S}$. The scheme ${\cal S}$ is canonically
isomorphic to $\ul{\rm Isom}(V\otimes\O_{\cal U^\p},\A_{C|{\cal U}^\p}^\vee)/\widetilde{G}$, i.e.,
$\sigma$ defines
$$
\tau^\p\colon {\Sym}^*(V\otimes\A_{C|{\cal U}^\p})^{\widetilde{G}} \lra \O_{\cal U^\p}
$$
which gives
$$
\widetilde{\tau}_C\colon {\Sym}^*(V\otimes\A_{C})^{\widetilde{G}}
\subset j_*\Bigl({\Sym}^*(V\otimes\A_{C|{\cal U}^\p})^{\widetilde{G}}\Bigr) 
\stackrel{j_*(\tau^\p)}{\lra} j_*(\O_{\cal U^\p})=\O_{C\times X}.
$$
By Lemma~\ref{honest}, the singular principal bundle $(\A_{C|\{*\}\times X}, \widetilde{\tau}_{C|\{*\}\times X})$
must be an honest one. Therefore, $w$ is the image of the corresponding point in $\frak Y^{\rm h}(\widetilde{G})$,
and this implies the claim. 
\par
Next, let $(\A_{\frak I},\phi_{\frak I})$ be the restriction of the universal family on 
${\frak W}_{\ul{a},\ul{b},\ul{c}}\times X$ to ${\frak I}\times X$. As usual, we let
${\cal U}_{\rm lf}$ be the open subset of $\frak I\times X$ where $\A_{\frak I}$ is locally free.
As above, $\A^\vee_{\frak I|{\cal U}_{\rm lf}}$ and $\GL(V)/\widetilde{G}\subset W$ define a subscheme
${\cal S}\subset \A^\vee_{\frak I|{\cal U}_{\rm lf}:\ul{a},\ul{b},\ul{c}}$, and the section
$\sigma^\p\colon {\cal U}_{\rm lf}\lra  \A^\vee_{\frak I|{\cal U}_{\rm lf}:\ul{a},\ul{b},\ul{c}}$ defined by
$\phi_{\frak I| \A^\vee_{\frak I|{\cal U}_{\rm lf}:\ul{a},\ul{b},\ul{c}}}$ factorizes over ${\cal S}$, so that
we find, as before, a section
$$
\sigma\colon  {\cal U}_{\rm lf}\lra \ul{\rm Isom}(V\otimes\O_{{\cal U}_{\rm lf}},\A^\vee_{\frak I|{\cal U}_{\rm lf}})
/\widetilde{G}
$$
which leads to a homomorphism
$$
\widetilde{\tau}_{\frak I}\colon {\Sym}^*(V\otimes\A_{\frak I})^{\widetilde{G}}
\lra \O_{\frak I\times X}.
$$
The family $(\A_{\frak I}, \widetilde{\tau}_{\frak Y})$ defines a $\GL(Y)$-equivariant morphism
$\frak u\colon \frak I\lra \frak t^{-1}\bigl({\frak W}^{\eps\rm-ss}_{\ul{a},\ul{b},\ul{c}}\bigr)$
with $\frak t^\p\circ \frak u=\id_{\frak I}$, and this is even stronger than the assertion of the theorem.
\end{proof}
\begin{Rem}
The construction outlined above may be used to replace the construction given in \cite{Schmitt}.
However, the construction in \cite{Schmitt} works also on a wide class of singular varieties, as
shown in \cite{Bhosle}, when ignoring the condition ``$\det(\A)\cong\O_X$".
\end{Rem}

\end{document}